\documentclass[11pt,letterpaper]{amsart}

\usepackage{amssymb}
\usepackage{amsthm}
\usepackage{amsfonts}
\usepackage{amsmath}
\usepackage{amstext}

\numberwithin{equation}{section}

\theoremstyle{plain}

\newtheorem{thm}{Theorem}[section]
\newtheorem{prop}[thm]{Proposition}
\newtheorem{cor}[thm]{Corollary}
\newtheorem{lem}[thm]{Lemma}

\theoremstyle{definition}

\theoremstyle{remark}
\newtheorem{rem}[thm]{Remark}

\newcommand{\one}{\mathbf {1}}
\newcommand{\algc}{\operatorname{A}}
\newcommand{\sign}{\mbox{$\operatorname {sgn}$}\,}
\DeclareMathOperator{\hlink}{\tilde \Lambda }
\DeclareMathOperator{\link}{ \Lambda }
\DeclareMathOperator{\alink}{\Omega }
\DeclareMathOperator{\ahlink} {\tilde \Omega }
\DeclareMathOperator{\supp}{\operatorname{supp} }

\newcommand{\R}{\mathbb{R}}

\newcommand{\Z}{\mathbb{Z}}
\newcommand{\Q}{\mathbb{Q}}

\newcommand{\inv}{^{-1}}

\newcommand{\cn}{\mathbf {a}}

\newcommand{\half}{\tfrac12}
\def\({(\!(}
\def\){)\!)}

\newcommand{\dchi}{\, d\chi}
\newcommand{\pol}{\mathcal{P}}
\newcommand{\pola}{\mathcal{A}}
\newcommand{\ideal}{\mathcal{I}}
\newcommand{\lpol}{\widetilde{\mathcal{P}}}

\date{September 7, 1998}

\thanks{Research supported by CNRS. First author also supported by NSF grant
DMS-9628522.}

\title[Real Algebraic Sets of Dimension 4]
{Topology of Real Algebraic Sets\\ of Dimension 4: \\
Necessary Conditions}

\author{Clint McCrory and Adam Parusi\'nski}

\address{D\' epartement de Math\' ematiques, Universit\' e d'Angers,
   2, bd Lavoisier, 49045 Angers cedex 01, France}
\email{parus@tonton.univ-angers.fr}
\address {Department of Mathematics, University of Georgia, Athens, GA
30602, USA }

\email{clint@math.uga.edu}


\newcommand{\abstracttext}{}

\begin{document}
\begin{abstract} \abstracttext
{Operators on the ring of algebraically constructible functions
are used to compute local obstructions for a four-dimensional
semialgebraic set to be homeomorphic to a real algebraic set. The
link operator and arithmetic operators yield $2^{43}-43$
independent characteristic numbers mod 2, which generalize the
Akbulut-King numbers in dimension three.}
 \end{abstract} \maketitle

The ring of algebraically constructible functions
on a real algebraic set was introduced in \cite{mccpar}. The link
operator on this ring was used to give a new description of the
Akbulut-King numbers of three-dimensional stratified sets
\cite{akbking2}, as well as to generalize the topological
conditions on algebraic sets discovered by Coste and Kurdyka
\cite{coskur}.

Akbulut and Kurdyka have asked what invariants can be constructed
by our method for four-dimensional sets.  In this paper we produce a large
number of independent new local topological conditions satisfied by
algebraic sets of dimension four.
Thus, in particular, there are four-dimensional semialgebraic sets which
have vanishing
Akbulut-King invariants, but which are not homeomorphic to algebraic
sets. We do not know whether a four-dimensional semialgebraic set
satisfying all of our conditions, as well as the Akbulut-King conditions,
must be homeomorphic to an algebraic set.

The properties of the link operator $\link$ on constructible functions
are reviewed in section \ref{Acf}. The main result of
\cite{mccpar} (see also \cite{parszaf1}, \cite{parszaf2}) is that
the operator $\hlink = \half\link$ preserves the set of algebraically
constructible functions: If $\varphi$ is algebraically constructible,
so is $\hlink\varphi$. Thus, in particular, if $X$ is a real
algebraic set with
characteristic function $\one_X$, every function obtained from
$\one_X$ using the arithmetic operations $+$, $-$, $*$, and the
operator $\hlink$ is integer-valued. Sets with this property we call
{\it completely euler}. The property that $\hlink\one_X$ is
integer-valued is equivalent to Sullivan's condition \cite{sul} that
$X$ is {\it euler}: for
all $x\in X$, the link of $x$ in $X$ has even euler
characteristic.

In section \ref{Cons} we show how to construct systematically invariants
of a set $X$ which vanish if and only if $X$ is completely euler.
These invariants are local, and if we assume that all the links
of points of $X$ are completely euler, then $X$ is completely
euler if and only if a finite list of mod 2 characteristic
numbers vanish for each link. For $X$ of dimension at most 4 we find
$2^{29}-29$ such characteristic numbers. Section \ref{Ind} contains a
proof of the independence of these numbers, by construction of
examples which distinguish them.

There are also arithmetic operators with rational coefficients
which preserve the set of algebraically constructible functions.
For example, if $\varphi$ is algebraically constructible, then so
is $P(\varphi)=\half(\varphi^4-\varphi^2)$. In section
\ref{Arith} we characterize such operators, and we show that $P$
is the only such operator which gives new conditions on the
topology of 4-dimensional algebraic sets. Using the operator $P$
together with the operator $\hlink$, we enhance our previous
construction to produce $2^{43}-43$ independent local characteristic
numbers which vanish for algebraic sets of dimension $\leq 4$.

We work with real semialgebraic sets in euclidean space,
semialgebraic Whitney stratifications, and semialgebraically
constructible functions. We use the foundational results that the
link of a point in a semialgebraic set is well-defined up to
semialgebraic homeomorphism, and that a
semialgebraic set has a semialgebraic, locally semialgebraically
trivial Whitney stratification ({\it
cf.~}\cite{coskur}). Since semialgebraic sets are triangulable and our
constructions are purely topological, we could just as well work with
piecewise-linear sets, stratifications, and constructible functions.

It is not natural to assume that a real semialgebraic set
has the same local dimension at every point, so we must be careful
in dealing with dimension. By the dimension of a semialgebraic
set $X$ we mean its topological dimension, denoted $\dim X$. If
$T$ is a semialgebraic triangulation of $X$, then $\dim X$ is the
maximum dimension of a simplex of $T$. For $x\in X$, the local
dimension $\dim_xX$ of $X$ at $x$ is the maximum dimension of a
closed simplex of $T$ containing $x$. Thus $\dim
X=\operatorname{max}\{\dim_xX,\, x\in X\}$, and if $L$ is the link
of $x$ in $X$, then $\dim L=\dim_xX-1$, which may be less than
$\dim X - 1$.

\medskip
\section{Algebraically Constructible Functions}\label{Acf}

\subsection{Definition and main properties}\label{Acf1}
Let $X\subset \R ^n$ be a real algebraic set. Following
\cite{mccpar} we say that an integer-valued function $\varphi: X\to \Z$ is
{\it algebraically constructible} if there exists a finite collection
of algebraic sets $Z_i$ and proper regular morphisms ${f_i}:Z_i \to X$
such that
$\varphi$ admits a presentation as a finite sum
\begin{equation}\label{e:algcons}
\varphi (x) = \sum m_i \chi (f_i\inv (x))
\end{equation}
with integer coefficients $m_i$, where $\chi$ is the euler
characteristic.   We recall from \cite
{parszaf1}, \cite{parszaf2} that  $\varphi: X\to \Z$ is algebraically
constructible if and only if there exists a finite set of
polynomials $g_1, \ldots, g_s \in \R[x_1, \ldots ,x_n]$ such that
\begin{equation}\label{e:signs}
\varphi (x) = \sign g_1(x) +  \ldots + \sign g_s(x) .
\end{equation}
The set of algebraically constructible functions on $X$ forms a ring,
which we
denote by $\algc (X)$.  An algebraically constructible function is
constructible in the usual sense; that is, it is a finite
combination of characteristic functions of semialgebraic subsets
$X_i$ of $X$,
\begin{equation} \label{e:cons}
\varphi = \sum n_i \one_{X_i}
\end{equation}
with integer coefficients $n_i$. The converse is not true; there
exist constructible functions which are not algebraically
constructible (see \cite{mccpar}).

Let $X$ be a semialgebraic set, and let $\varphi$ be a constructible
function on  $X$ given by
(\ref{e:cons}).  Without loss of generality we may assume that all
$X_i$ are closed in $X$.  If, moreover, all
$X_i$ are compact then we define the {\it euler integral} of $\varphi$ by
\begin{equation*}
\int \varphi \dchi = \sum n_i \chi(X_i).
\end{equation*}
It is easy to see that $\int \varphi \dchi$ is well-defined; it
depends only on $\varphi$ and not on the presentation (\ref{e:cons}).
The euler integral is uniquely determined by the properties that
it is linear,
\begin{equation*}
\int \varphi + \psi \dchi = \int \varphi\dchi + \int\psi \dchi,
\end{equation*}
and the euler integral of the characteristic function of a compact
semialgebraic set is its euler characteristic,
\begin{equation*}
\int \one_W \dchi = \chi(W).
\end{equation*}
If $\varphi$ is given by (\ref{e:cons}), then for $x\in X$,
the {\it link} of $\varphi$ at $x$ is defined by
\begin{equation*}
\link \varphi (x) = \sum n_i \chi(X_i \cap S(x,\varepsilon)),
\end{equation*}
where $ S(x,\varepsilon)$ is a sufficiently small sphere centered at $x$.
Thus
\begin{equation*}
\link \varphi (x) = \int_L\varphi\dchi,
\end{equation*}
where $L=X \cap S(x,\varepsilon)$ is the link of $x$ in $X$.
It is
easy to see that $\link \varphi$ is a well-defined constructible
function on $X$. The link operator $\link$ is linear, and it
satisfies
\begin{equation} \label{e:ll}
\link\link = 2\link
\end{equation}
(see \cite{mccpar}).
If $\varphi$ has compact support then the euler integral of the link of
$\varphi$ is zero \cite{mccpar},
\begin{equation} \label{e:intl}
\int \link\varphi\dchi = 0.
\end{equation}

\subsection{Completely euler spaces}
Let $\hlink = \half \link$. The following result is proved in
\cite{mccpar}; see also \cite{parszaf1}.

\begin{thm}\label{t:hlink}
If $\varphi \in \algc (X)$, then
 $\link \varphi$ is an even-valued function and
$\hlink \varphi   \in \algc (X)$.
\end{thm}

Suppose now that the semialgebraic set $X$ is  homeomorphic to
an algebraic set.  Then by Theorem \ref{t:hlink} all
functions constructed from $\one_X$ by means of the arithmetic
operations $+,-,*$, and the operator $\hlink$, are integer-valued.
This observation was used in \cite{mccpar} to recover the Akbulut-King
conditions \cite{akbking2} on the topology of real algebraic sets
of dimension $\leq 3$.

We say that the constructible function $\varphi$ on $X$ is {\it euler}
if $\link\varphi(x)$ is even for all $x\in X$, {\it i.e.}
$\hlink\varphi$ is integer-valued. The semialgebraic set
$X$ is {\it euler} if $\one_X$ is euler. We say that $\varphi$ is
{\it completely euler} if all functions obtained from $\varphi$ by
using the operations $+$, $-$, $*$, $\hlink$ are integer-valued.
The semialgebraic set
$X$ is {\it completely euler} if $\one_X$ is completely euler. Thus
a necessary condition for a semialgebraic set $X$ to be homeomorphic to
an algebraic set is that $X$ is completely euler.

The following theorem follows easily from the theory of basic open sets
in real algebraic geometry \cite{bcr};  for the proof see \cite {mccpar}.

\begin{thm}\label{t:2d}
If $X$ is an algebraic set of dimension at most $d$, then every
constructible function on $X$ divisible by $2^d$ is algebraically
constructible.
\end{thm}

The following related observation will be useful in the sequel.
If $X$ is a semialgebraic set, let $\ideal (X)$ be the set of
constructible functions $\varphi$ on $X$ such that, for all $k\geq 0$,
$\dim(\supp \varphi \pmod {2^k}) < k$. Note that $\ideal (X)$ is an
ideal of the ring of constructible functions on $X$.

\begin{prop}\label{p:ideal}
$\hlink\ideal (X)\subset \ideal (X)$.
\end{prop}

The proof is elementary ({\it cf}.~ \cite{mccpar}, Prop. 3.1). It is easy
to see that Theorem \ref{t:2d}
implies that $\ideal (X)\subset \algc(X)$, while Theorem
\ref{t:hlink} says that $\hlink \algc(X)\subset \algc(X)$.

\subsection{Properties of $\hlink$ and $\ahlink$}\label{properties}
Define an operator $\alink$ on constructible functions by
$\alink\varphi= 2\varphi - \link\varphi$, and set $\ahlink =
\half\alink$. Thus
\begin{equation*}
\hlink + \ahlink = \operatorname I,
\end{equation*}
where $\operatorname I$ is the identity operator.
The following
elementary properties of $\hlink$ and $\ahlink$ will be used
repeatedly:

(a) {\bf Composition}:
\begin{equation*}
\hlink\hlink = \hlink,\ \ \ahlink\ahlink = \ahlink,\ \ \hlink \ahlink=0,
\ \ \ahlink\hlink = 0.
\end{equation*}

(b) {\bf Integral}:
\begin{equation*}
\int\hlink\varphi\dchi = 0,\ \ \int\ahlink\varphi\dchi = \int\varphi\dchi.
\end{equation*}

(c) {\bf Support}:

Let $l$ be a nonnegative integer.
If $\dim(\supp\varphi)\leq 2l$, then $\dim(\supp(\hlink\varphi))
\leq 2l - 1$.
If $\dim(\supp\varphi)\leq 2l+1$, then
$\dim(\supp(\ahlink\varphi))\leq 2l$.

\vskip.1in
(d) {\bf Slice}:

Let $W\subset X$ be given by $W=H\cap X$, where $X$ is a semialgebraic
set in $\R^n$ and $H$
is a smooth $(n-1)$-dimensional semialgebraic subset of $\R^n$ which is
transverse to a semialgebraic Whitney stratification $\mathcal S$ of $X$.
If $\varphi$
is a constructible function on $X$ which is constant on the
strata of $\mathcal S$, then
\begin{equation*}
(\hlink\varphi)|W = \ahlink(\varphi|W),\ \
(\ahlink\varphi)|W = \hlink(\varphi|W).
\end{equation*}

Properties (a) and (b) follow from the corresponding properties of
$\link$ ((\ref{e:ll}), (\ref{e:intl})). Properties (c) and (d) follow from
the definition of $\link$ by the following topological arguments.

Given a constructible function $\varphi$ on $X$, let
$\mathcal S$ be a locally trivial stratification of $X$ such that
$\varphi$ is constant on the strata of $\mathcal S$. Let
$k =\dim(\supp\varphi)$, and let $x\in S$, where $S$ is a $k$-dimensional
stratum of $\mathcal S$. Then
\begin{equation*}
\hlink\varphi(x) = \half\int_{X\cap S(x,\epsilon)}\varphi\dchi
=\half\int_{S\cap S(x,\epsilon)}\varphi\dchi.
\end{equation*}
For $\epsilon >0$ sufficiently small,
$S\cap S(x,\epsilon)$ is semialgebraically homeomorphic to a $(k-1)$-sphere.
Thus if $k$ is even, $\hlink\varphi(x)=0$; and if
$k$ is odd, $\hlink\varphi(x)=\varphi(x)$, so $\ahlink\varphi(x)=0$.

Since restriction and the operators $\hlink$ and $\ahlink$ are
linear, by the definition of $\link\varphi$ it suffices to prove
(d) for $\varphi=\one_X$. For all $x\in W$, $X\cap S(x,\epsilon)$ is
semialgebraically
homeomorphic to the suspension of $W\cap S(x, \epsilon)$,
since a neighborhood of $x$ in $X$ is
semialgebraically homeomorphic to the product of
an interval with a neighborhood of $x$ in $W$. Thus
\begin{equation*}
\begin{split}
\chi(X\cap S(x,\epsilon))&=2 - \chi(W\cap S(x,\epsilon)),\\
\link\one_X(x)&=2\cdot\one_W(x)-\link\one_W(x),\\
\link\varphi(x)& = 2(\varphi|W)(x) - \link(\varphi|W)(x),
\end{split}
\end{equation*}
which implies (d).


\medskip.
\section{Construction of local invariants from the link operator}
\label{Cons}

\subsection{A general construction}\label{Cons1}

Let $X\subset \R^n$ be a semialgebraic set. We describe an algorithm for
producing
a list of $\Z/2$-valued obstructions (finite in number if $X$ has
a finite locally trivial stratification) the
vanishing of which are
necessary and sufficient for $X$ to be completely euler. Let
$\mathcal S$ be a locally trivial stratification of $X$, and let
$X_0\subset X_1\subset X_2 \subset\cdots\subset X_d=X$, $d=\dim X$,
be the skeletons of $\mathcal S$. Let $\hlink(X)$ be the ring of
all functions obtained from $\one_X$ by the operations $+$, $-$,
$*$, $\hlink$.
Every such function is constant on the strata of $\mathcal S$ and
takes rational values. By definition the set $X$ is completely euler if and
only if all functions in $\hlink(X)$ are integer-valued.

We define a sequence of subrings of $\hlink(X)$,
\begin{equation*}
\hlink_0(X)\subset\hlink_1(X)\subset\hlink_2(X)\subset\cdots,
\end{equation*}
inductively as follows. Let $\hlink_0(X)$ be the ring generated by
$\one_X$. For $k\geq0$ let
$\hlink_{k+1}(X)$ be the ring generated by $\hlink_k(X)$ and
$\{\hlink\varphi\ |\ \varphi\in\hlink_k(X)\}$.  Then $\hlink(X)$ is
the direct limit of the rings $\hlink_k(X)$. For each $k\geq 0$, if
all functions in $\hlink_k(X)$ are integer-valued, we define
obstructions for all functions in $\hlink_{k+1}(X)$ to be integer-valued,
and we show that if all functions in $\hlink_d(X)$ are
integer-valued, then so are all functions in $\hlink(X)$.

Recall that $\ideal (X)$ is the ideal of constructible functions
$\varphi$ on $X$ such that, for all $k\geq 0$,
$\dim(\supp \varphi \pmod {2^k}) < k$.
Let $k\geq 0$ and suppose that $\hlink_k(X)$ is integer-valued.
Moreover, suppose we have finite subsets $G_0,\dots,G_{k-1}$ of
$\hlink_k(X)$ such that for $j=0,\dots,k-1$, $G_0\cup
\cdots\cup G_j$ is an additive generating set for
$\hlink_j(X)/(\ideal (X)\cap\hlink_j(X))$, and for all
$\varphi\in G_j$, $\supp\varphi\subset X_{d-j}$. Then
$\hlink_k(X)/(\ideal (X)\cap\hlink_k(X))$ is additively
generated by $G_0, \dots,G_{k-1}$, together with the products
\begin{equation}\label{e:phi}
\Phi = \begin{cases}
\varphi\prod(\hlink\varphi_i)^{n_i},& d-(k-1)\ \text{even},\\
\varphi\prod(\ahlink\varphi_i)^{n_i},& d-(k-1)\ \text{odd},
\end{cases}
\end{equation}
where $\varphi\in G_0\cup\cdots\cup G_{k-1}$, $\varphi_i\in
G_{k-1}$, and $n_i>0$ for some $i$.
By the support properties of $\hlink$ and $\ahlink$
(section \ref{properties}(c)),
$\supp\Phi\subset X_{d-k}$, so by definition of $\ideal (X)$
the quotient $\hlink_k(X)/(\ideal (X)\cap\hlink_k(X))$ is
additively generated by $G_0\cup\cdots\cup G_{k-1}\cup G_k$,
where $G_k$ is a finite set of representatives of the equivalence
classes modulo $2^{d-k}$ of the functions $\Phi$.

Recall that $\hlink(\ideal (X))\subset\ideal (X)$. Thus if
 all functions in $\hlink_k(X)$ are integer-valued, then all
functions in $\hlink_{k+1}(X)$ are integer-valued if and only if all the
functions $\varphi\in G_k$ are euler. We call the conditions that
these functions are euler the {\it depth} $k+1$ {\it euler conditions}
for $X$.
Since the construction of the generating sets $G_k$ terminates at
$k=d$ (for $k=d+1$, $\supp\Phi=\emptyset$), we have that
$\hlink_d(X)/(\ideal (X)\cap\hlink_d(X)) =
 \hlink(X)/(\ideal (X)\cap\hlink(X))$.
Thus $X$ is completely euler if and only if the depth $k$ euler
conditions hold for all $k=1,\dots,d$.

\subsection{Localization of the invariants}\label{Cons2}
The condition that $X$ is completely euler is local. The
constructible
function $\varphi$ on $X$ is euler if and only if
$\int_L\varphi\dchi\equiv
0 \pmod 2$ for all links $L$ in $X$, so each euler condition for $X$
is
equivalent to the vanishing of a set of mod 2 obstructions. (If $X$
has a finite locally trivial
stratification then the number of semialgebraic
homeomorphism types of links in $X$ is finite,
so the total number of obstructions is finite.)

It is easier to analyze these obstructions if we restrict our attention
to a single link $L$ in $X$. This has the advantage of reducing the
dimension of the space we consider.
By the slice property of the link operator (section \ref{properties}(d)),
 if the constructible function
$\varphi$ on $X$ is euler, then $\varphi |L$ is euler and $\int_L\varphi
\dchi$ is even.  Thus the compact semialgebraic set $L$ is the link of
a point in a completely euler space $X$ if and only if the following
conditions are satisfied:

(A) $L$ is completely euler, {\it i.e.} for all $\varphi\in\hlink(L)$,
$\varphi$ is integer-valued.

(B) For all $\varphi\in\hlink(L)$, $\int\varphi\dchi$ is even.

We have seen that (A) is equivalent to a finite list of euler conditions
for $L$. In the
above discussion $X$ is replaced by $L$, and $d=\dim X$ is
replaced by $d'=\dim L$ (with $d'\leq d-1$ if $L$ is a link in $X$).
 If we assume that $L$ is completely euler,
then condition (B) can be reduced to the vanishing of a finite set
of mod 2 invariants. This follows from the facts that the ring
$\hlink(L)/(2\ideal (L)\cap\hlink(L))$ is additively finitely
generated, and that if $\varphi\in 2\ideal (L)$ then $\int\varphi\dchi$
is even. (Note that if $L$ is a link in $X$ then
 $2\ideal (L)$ is the restriction to $L$ of the
ideal $\ideal (X)$.) To produce a finite set of additive generators
for this ring, note as above that
if $k\geq 1$ and we have finite subsets $G_0,\dots,G_{k-1}$ of
$\hlink(L)$ such that for $j=0,\dots,k-1$, $G_0\cup\cdots\cup
G_j$ is an additive generating set for
 $\hlink_j(X)/(2\ideal (X)\cap\hlink_j(X))$ and for all
$\varphi\in G_j$ $\supp\varphi\subset L_{d'-j}$, then
 $\hlink_k(X)/(2\ideal (X)\cap\hlink_k(X))$ is additively
generated by $G_0\cup\cdots\cup G_{k-1}\cup G_k$, where $G_k$ is
a finite set of representatives of the equivalence classes modulo
$2^{d'-k+1}$ of the functions $\Phi$ \eqref{e:phi}.
We call the integrals
$\int\psi\dchi \pmod 2$, $\varphi\in G_k$ the {\it depth k characteristic
numbers} of
$L$. If $L$ is completely euler, then $\int\varphi\dchi$ is even for
all $\varphi\in\hlink(L)$ if and only if the depth $k$ characteristic
numbers of $L$ vanish for all $k=0,\dots,d'-1$.

\subsection{Invariants in dimension 4}\label{Inv4}
Now we work out the euler conditions and characteristic numbers for
a link $L$ in a semialgebraic set $X$ of dimension at most 4.

(A) {\it $L$ is completely euler.} The computation of the euler conditions
for a semialgebraic set of dimension at most 3 was done in \cite{mccpar}.
We review it here as an illustration of the above general algorithm.
The ring $\hlink_0(L)$ is generated by $\one_L$, and the depth 1
euler condition is that $\one_L$ is euler (Sullivan's condition).
 Assume this, and let
$\varphi = \ahlink\one_L\in\hlink_1(L)$. Then
$\hlink_1(L)/(\ideal (L)\cap\hlink_1(L))$ is additively generated
by $\one_L$, $\varphi$, $\varphi^2$, $\varphi^3$, since
$\supp \varphi\subset
L_2$ and $\varphi^4\equiv \varphi^2\pmod 4$.  So the depth 2 euler
conditions are that $\varphi$, $\varphi^2$, $\varphi^3$ are euler. But
$\varphi$ is euler since $\hlink\varphi =\hlink\ahlink\one_L=0$,
and $\varphi^2$ and $\varphi^3$ are also euler since
$\varphi^2$ and $\varphi^3$ are congruent to $\varphi$ mod 2.
Thus the depth 2 euler conditions are trivial.

Let $\beta =\hlink(\varphi^2)$ and $\gamma=\hlink(\varphi^3)$. The
ring $\hlink_2(L)/(\ideal (L)\cap\hlink_2(L))$ is additively
generated by $\one_L$, $\varphi$, $\varphi^2$, $\varphi^3$, and the
products
$\psi_{abc}= \varphi^a\beta^b\gamma^c$ with $b+c>0$. Since
$\supp(\psi_{abc})\subset L_1$, we may write
\begin{equation*}
\psi_{abc}=\alpha^a\beta^b\gamma^c
\end{equation*}
where $\alpha = \varphi |L_1$. Thus
$\hlink_2(L)/(\ideal (L)\cap\hlink_2(L))$
is additively generated by $\one_L$, $\varphi$, $\varphi^2$, $\varphi^3$,
together with representatives of the equivalence classes mod 2 of
the functions $\psi_{abc}$, namely $\beta$, $\gamma$, $\alpha\beta$,
$\alpha\gamma$,
$\beta\gamma$, $\alpha\beta\gamma$. Now $\beta$ and $\gamma$ are euler,
since $\hlink\beta=\hlink\hlink(\varphi^2)=\hlink(\varphi^2)=\beta$,
and similarly $\hlink\gamma=\gamma$. So the depth 3 euler conditions
on $L$ are that the four functions
$\alpha\beta$, $\alpha\gamma$, $\beta\gamma$, $\alpha\beta\gamma$,
are euler. (It was proved in \cite{mccpar} that these four functions
are euler for a semialgebraic set of dimension $\leq 3$ if and only if the
Akbulut-King conditions hold.) We conclude that $L$ is completely
euler if and only if the functions
\begin{equation*}
\one_L,\ \varphi\beta,\ \varphi\gamma,\ \beta\gamma,\ \varphi\beta\gamma
\end{equation*}
are euler, where $\varphi=\ahlink\one_L$, $\beta=\hlink(\varphi^2)$,
and $\gamma=\hlink(\varphi^3)$.

(B) {\it For all $\varphi\in\hlink(L)$, $\int\varphi\dchi$ is even.} There
is one depth 0 characteristic number, $\int\one_L\dchi \pmod 2$.
(The vanishing of this characteristic number is Sullivan's condition
$\chi(L)\equiv 0 \pmod 2$ for a link in an algebraic set.) Again consider
$\varphi =\ahlink\one_L$. The ring
$\hlink_1(L)/(2\ideal (L)\cap\hlink_1(L))$
is additively generated by
$\one_L$, $\varphi$, $\varphi^2$, $\varphi^3$, since $\supp \varphi\subset
L_2$ and $\varphi^4\equiv \varphi^2+ 2(\varphi-\varphi^3)\pmod 8$.
So the depth 1 characteristic numbers are the euler integrals
of $\varphi$, $\varphi^2$, and $\varphi^3 \pmod 2$. But $\int\varphi
\dchi=\int\ahlink\one_L\dchi=\int\one_L\dchi$, and $\varphi^2$,
$\varphi^3$ have the same euler integrals  as $\varphi$ mod 2, since they
are
congruent to $\varphi$ mod 2. So there are no new characteristic
numbers of depth one.

Now consider $\beta=\hlink(\varphi^2)$ and $\gamma=\hlink(\varphi^3)$.
The ring $\hlink_2(L)/(2\ideal (L)\cap\hlink_2(L))$ is generated
additively by $\one_L$, $\varphi$, $\varphi^2$, $\varphi^3$, together with
representatives of the equivalence classes mod 4 of the products
$\psi_{abc} = \alpha^a\beta^b\gamma^c$ with $b+c>0$. (Here $\alpha =
\varphi |L_1$ as before.) So it suffices to consider the exponents
$a, b, c$ equal to $0, 1, 2, 3$. This still leaves us with $60 =
4^3-4$ functions.
Note that $\beta \equiv \beta ^2 \equiv \beta ^3 \pmod 2$, so instead of
$\beta$, $\beta ^2$, $\beta ^3$ we use
$\beta$, $\beta_2=\beta ^2 - \beta$, $\beta_3= \beta ^3- \beta$ as a
set of additive generators for the ring generated by $\beta$ modulo $4$.
The advantage of this set of generators
is that $\beta_2$ and $\beta_3$ are even-valued.  We let $\beta_1=\beta$.
Similarly we define
$\alpha_1$, $\alpha_2$, $\alpha_3$, $\gamma_1$, $\gamma_2$, $\gamma_3$.
Now
all products $\psi'_{abc}=\alpha_a\beta_b\gamma_c$ with two even-valued
factors are divisible by $4$.
Therefore the following list of functions completes our additive
generating set:

\begin{gather*}
 \beta,\, \gamma; \tag{$S_1$}\\
 \alpha \beta ,\, \alpha\gamma ,\, \beta \gamma, \alpha \beta \gamma;
\tag{$S_2$}\\
 \beta_2 ,\, \gamma_2,\, \beta_3,\, \,\gamma_3 ;
\tag{$S_3$}
\end{gather*}

\begin{gather*}\tag{$S_4$}
\alpha \beta_2 ,\, \alpha \gamma_2,\, \beta \alpha_2, \, \beta
\gamma_2,\, \gamma\alpha_2, \, \gamma\beta_2, \\
 \alpha\beta_3, \, \alpha\gamma_3, \, \beta\alpha_3, \,
\beta\gamma_3, \, \gamma\alpha_3, \, \gamma\beta_3,\\
 \alpha\beta\gamma_2, \, \alpha\gamma\beta_2, \, \beta\gamma\alpha_2, \\
  \alpha\beta\gamma_3,\, \alpha\gamma\beta_3, \, \beta\gamma\alpha_3.
\end{gather*}

The characteristic numbers of depth 2 are the euler integrals mod 2
of these functions.
Now $\int\beta\dchi = \int\hlink(\varphi^2)\dchi = 0$ and
$\int\gamma\dchi = \int\hlink(\varphi^3)\dchi = 0$, so the
integrals of the functions in $S_1$ are zero.
The functions in $S_3$ or $S_4$ are even, and hence have even euler
integrals. The
nontrivial characteristic numbers of depth 2 are the integrals of the
functions in $S_2 \pmod 2$:
\begin{equation}\label{akint}
\int\alpha \beta\dchi ,\, \int\alpha\gamma\dchi ,\,
\int\beta \gamma\dchi,\, \int \alpha \beta \gamma\dchi
\end{equation}
The vanishing of these 4 characteristic numbers is a natural
generalization of the Akbulut-King conditions.

Let $S = S_1\cup S_2\cup S_3\cup S_4$, and let $S'= S_2\cup S_3\cup S_4$.
The ring $\hlink_3(L)/(2\ideal (L)\cap\hlink_3(L))$
is generated additively by $\one_L$, $\varphi$, $\varphi^2$, $\varphi^3$,
and $S$, together with the equivalence classes mod 2 of the functions
\begin{equation}
  \label{e:prod}
\Phi(\mathbf m, \mathbf n)=
  \alpha ^{m_\alpha}\beta ^{m_\beta}\gamma ^{m_\gamma}
\prod_{\psi\in S'} (\ahlink (\psi)) ^{n_\psi} ,
\end{equation}
where $\mathbf m = (m_\alpha,m_\beta,m_\gamma)$, and $\mathbf n =
(n_\psi)_{\psi\in S'}$, with $m_\alpha$, $m_\beta$, $m_\gamma$,
$n_\psi$ equal to 0 or 1, and $\sum_{\psi\in S'}n_\psi>0$.
 The support of
$\Phi(\mathbf m, \mathbf n)$ is contained in $L_0$.
Since  $L_0$  is finite,
\begin{equation} \label{e:inv}
  \int \Phi(\mathbf m, \mathbf n) \dchi =
\sum_{p\in L_0} \Phi(\mathbf m, \mathbf n)(p) .
\end{equation}
Some of these integrals are automatically even.
For instance,
\begin{eqnarray*}
  \int \ahlink \beta_2 \dchi = \int \beta_2 \dchi \equiv 0 \pmod 2,
\end{eqnarray*}
 since $\beta_2$ is even-valued.  The same argument works for any
other function from $S_3$ or $ S_4$. On the other hand, for
$\psi \in S_2$ we recover the depth 2 characteristic numbers
$\int \ahlink  \psi\dchi  = \int \psi \dchi$ (\ref{akint}). This leaves
us with $2^3(2^{26}-1) - 26= 2^{29}-34$ new characteristic numbers
(\ref{e:inv}) of depth 3.  We summarize our analysis of the link $L$
as follows.

\begin{thm}\label{t:inv1}
Let $L$ be a compact semialgebraic set of dimension at most $3$.
Let $\varphi = \ahlink \one_L$, $\beta = \hlink \varphi ^2$,
 $\gamma = \hlink \varphi ^3$ . Then $L$ is
a link in a completely euler space if and only if the following
conditions hold:
\begin{enumerate}
\item
$L$ is euler and  $\varphi \beta$, $\varphi \gamma$,
$\beta \gamma$,  $\varphi \beta \gamma$
 are euler ($L$ satisfies the Akbulut-King conditions);
\item
$L$ has even euler characteristic and $\varphi \beta$,
$\varphi \gamma$, $\beta \gamma$, $\varphi \beta \gamma$  
have even  euler integral;
\item
The characteristic numbers
\begin{eqnarray*}
\cn (\mathbf m, \mathbf n) =
  \sum_{p\in L_0} \Phi(\mathbf m, \mathbf n)(p),
\end{eqnarray*}
 defined above (\ref{e:inv}) are even.
\end{enumerate}
\end{thm}
The completely euler set
 $L$ of dimension at most 3 has $2^{29}-29$ characteristic numbers
 in all: the euler characteristic
 of $L$ and the $2^{29}-30$ characteristic numbers (\ref{e:inv}), which
 include the 4 numbers (\ref{akint}).

\medskip
\section{Independence of invariants.  Construction of examples}\label{Ind}

In this section we show that the $2^{29} - 29$ characteristic numbers of
Theorem \ref {t:inv1} are independent, {\it i.e.} for any subset $\mathcal
C$ of these characteristic numbers there exists
 a compact completely euler semialgebraic set $L$ of dimension
$3$  such that exactly those
characteristic numbers of $L$ which are in $\mathcal C$ are nonzero
(mod 2).
  Since the characteristic numbers are additive
with respect to disjoint union,  it suffices to construct for each
characteristic number
a set $L$ with exactly that characteristic number nonzero.

\subsection{Outline of the construction.}
The set $L=S^3\vee S^3$ (two 3-spheres identified at at point $p$)
has the characteristic number
$\chi(L)\equiv 1 \pmod 2$, and $\varphi = \ahlink\one_L = \one_{\{p\}}$.
Thus $\beta=\hlink\varphi^2 = 0$ and $\gamma=\hlink\varphi^3 = 0$,
so $L$ is completely euler and all the other characteristic numbers
of $L$ are zero.

So it remains to construct, for each characteristic number
(\ref{e:inv}), a compact semialgebraic set $L$ of dimension 3 such that
$L$ has even euler characteristic, $L$ is euler, the functions
$\varphi \beta$, $\varphi \gamma$, $\beta \gamma$, $\varphi \beta \gamma$
are euler, and exactly the given characteristic number of $L$ is nonzero.

  Such a set $L$
will be constructed together with its filtration by skeletons of
a stratification
\begin{equation*}
L_0\subset L_1\subset L_2 \subset L_3=L
\end{equation*}
in the following three steps. We make no attempt to minimize the
number of strata in the construction of $L$.

\begin {enumerate}
\item
Given a characteristic number $\cn (\mathbf m, \mathbf n)$,
we construct $L_1$ and constructible functions $\alpha, \beta,
\gamma$ on $L_1$ such that:
        \begin {enumerate}
        \item
        $\ahlink \alpha = \ahlink \beta =\ahlink \gamma = 0$
        and $\alpha \beta, \alpha
        \gamma, \beta \gamma$, $ \alpha \beta \gamma$, and $\one_{L_1}$
        are euler,
        \item
        $\sum_{p\in L_0} \Phi(\mathbf m, \mathbf n)(p)$
        is the only sum as in (\ref {e:inv}) which is nonzero (mod 2).
        \end{enumerate}
The conditions that $\ahlink \alpha = 0$ and $\one_{L_1}$ is euler are not
necessary, but they will help to simplify the construction of $L$.
\item
Given $(L_1, \alpha, \beta, \gamma)$ satisfying 1(a), we construct $L_2$
and a constructible function $\varphi$ on $L_2$ such that:
        \begin {enumerate}
        \item
        $\hlink \varphi =0$ ,
        \item
        $\alpha \equiv \varphi |_{L_1}$, $\beta \equiv \hlink \varphi^2$,
        $\gamma
        \equiv \hlink \varphi^3 \pmod 4$ on open 1-dimensional strata,
        and the same congruences hold (mod 2) on $L_0$ ,
        \item
        $\int \varphi \dchi$ is even .
        \end{enumerate}
\item
Given $(L_2, \varphi)$ satisfying 2(a), we construct $L$ such that:
        \begin {enumerate}
        \item
        $L$ is euler,
        \item
         $\varphi\equiv \ahlink \one _L \pmod 8$ on open 2-dimensional
strata, $\varphi\equiv \ahlink \one_L\pmod 4$ on open 1-dimensional
         strata, and $\varphi\equiv \ahlink \one_L \pmod 2$ on $L_0$.
         \end {enumerate}
Then, in particular, if $L_2$ also satisfies 2(c) then $\chi(L)$ is even:
 $\int \one _L \dchi = \int \ahlink \one _L \dchi \equiv \int
\varphi \dchi \pmod 2$.
\end {enumerate}

Note that if $\psi\equiv\psi'\pmod{2^k}$, then
$\hlink\psi\equiv\hlink\psi'\pmod { 2^{k-1}}$.  This implies that
 the characteristic number $\cn(\mathbf m,\mathbf n)$
depends only on the congruence classes of $\one_L \pmod {16}$,
$\varphi\pmod 8$, or $\alpha$, $\beta$, $\gamma\pmod 4$.
We may also conclude the sightly stronger natural condition that
$\cn (\mathbf m,\mathbf n)$ depends only on the
congruence classes of $\one_L$, $\varphi$, $\alpha$, $\beta$,
$\gamma$ modulo the ideal $2\ideal (L)$. Indeed, the conditions 2(b)
and 3(b) above are modulo $2\ideal (L)$.

\subsection{Step 1. Construction of $(L_1, \alpha, \beta, \gamma)$. }
\label{Step1}

As noted above, we will construct $(L_1, \alpha, \beta,
\gamma)$ with the extra condition $\ahlink \alpha = 0$.  This
will allow
us to treat $\alpha, \beta, \gamma$ in the same way.  All the
characteristic numbers are defined modulo $2$;
on the other hand,  the
integer-valued functions $\alpha, \beta, \gamma$ on $L_1$ matter
 modulo $4$.

\subsubsection {Elementary blocks, an example} \label{blockexample}
First, to illustrate the method, we consider only one function, say
$\beta$.    In this case there are 4 characteristic
numbers:
\begin{equation}\label{e:aknum}
\begin{split}
&  \sum_{p\in L_0} \beta (p)\ahlink \beta_2 (p), \quad
 \sum_{p\in L_0} \beta (p)\ahlink \beta_3 (p), \\
&  \sum_{p\in L_0} \ahlink \beta_2 (p)\ahlink \beta_3 (p), \quad
  \sum_{p\in L_0} \beta (p) \ahlink \beta_2 (p) \ahlink \beta_3 (p) ,
\end{split}
\end{equation}
where $\beta_2=\beta ^2 - \beta$, $\beta_3= \beta ^3- \beta \pmod 2$.

Let $(L_1, \beta)$ be such that $\ahlink \beta = 0$. The characteristic
 numbers (\ref{e:aknum}) can be computed as follows.
Let $p\in L_1$ be a
vertex and suppose $p$ is in the boundary of exactly $k$
one-dimensional strata $S_1, \ldots, S_k$, with the
values of $\beta$ on these strata equal to
 $l_1, \ldots,l_k$, respectively.  Then,
since $\ahlink \beta(p)= \beta (p)- \half \sum l_i = 0$, we have
\begin{equation}\label {e:betas}
\begin{split}
& \beta(p) = \half \sum l_i, \\
& \ahlink \beta_2(p) = \ahlink \beta^2(p) =
\beta ^2 (p) - \half \sum l^2_i, \\
& \ahlink \beta_3(p) =  \ahlink \beta^3(p) =
\beta ^3 (p) - \half \sum l^3_i .
\end{split}
\end{equation}

Now consider the following 1-dimensional stratified sets
(with distinguished vertex $p$), which we call {\it  elementary blocks}:

\unitlength = 1cm
\begin{picture}(4,3)
\thicklines
\put(2.5,1){\circle{4}}
\put(1.8,1){\circle*{.2}}
\put(3.2,1){\circle*{.2}}
\put(0,0.9) {\shortstack{($A$)}}
\put(2,1.9) {\shortstack{$\beta = 1$}}
\put(2,-.1) {\shortstack{$\beta = 1$}}
\put(3.4,1) {\shortstack{$p$}}
\put(6,1) {\shortstack{$(\beta, \ahlink
\beta_2, \ahlink \beta_3)(p) = (1,0,0)$}}
\end{picture}

\begin{picture}(4,3)
\thicklines
\put(2.5,1){\circle{4}}
\put(1.8,1){\circle*{.2}}
\put(3.2,1){\circle*{.2}}
\put(0,0.9) {\shortstack{($B$)}}
\put(2.4,1.9) {\shortstack{$1$}}
\put(2.3,-.1) {\shortstack{$-1$}}
\put(3.4,1) {\shortstack{$p$}}
\put(6,1) {\shortstack{$(\beta, \ahlink
\beta_2, \ahlink \beta_3)(p) = (0,1,0)$}}
\end{picture}

\begin{picture}(4,3)
\thicklines
\put(2.5,1){\circle{4}}
\put(1.8,1){\circle*{.2}}
\put(3.2,1){\circle*{.2}}
\put(0,0.9) {\shortstack{($C$)}}
\put(2.4,1.9) {\shortstack{$2$}}
\put(2.4,-.1) {\shortstack{$0$}}
\put(3.4,1) {\shortstack{$p$}}
\put(6,1) {\shortstack{$(\beta, \ahlink
\beta_2, \ahlink \beta_3)(p) = (1,1,1)$}}
\end{picture}

\medskip

For each block we indicate the values of $\beta$ on the 1-dimensional
strata.  The values of $\beta$, $\ahlink \beta_2$, $\ahlink \beta_3$,
given by \eqref{e:betas}, are indicated to the right  of
each diagram.

\subsubsection {Additivity}
The numbers in \eqref{e:aknum} are additive with respect to taking the
wedge at the vertex $p$.  More precisely,
suppose we have two 1-dimensional stratified sets $(L', \beta')$
and $(L'', \beta'')$, with  $\beta'$ and $\beta''$  constructible
functions on $L'$ and $L''$ respectively, $\ahlink \beta' =0$, and
$\ahlink \beta'' =0$.  Suppose we have chosen vertices  $p'\in L'$ and
$p''\in L''$ as base points, and let $(L,p)$ be the wedge
$(L', p')\vee (L'',p'')$, with $p\in L$ the
induced base point on $L$.  Then $\beta'|_{L'\setminus p'}$ and
$\beta''|_{L''\setminus p''}$  extend uniquely to a
constructible function $\beta$  on $L$ satisfying $\ahlink \beta = 0$, and
\begin{equation}\label {e:add}
\begin{split}
& \beta (p) = \beta' (p') + \beta'' (p'') \\
& \ahlink \beta_2(p) \equiv \ahlink \beta'_2(p') +\ahlink \beta''_2(p'')
\pmod 2\\
& \ahlink \beta_3(p) \equiv  \ahlink \beta'_3(p') + \ahlink
\beta''_3(p'')
\pmod 2 .
\end{split}
\end{equation}

By wedging the elementary blocks ($A$), ($B$), ($C$), one may easily
construct
 $(L_1, \beta)$ with $\ahlink \beta = 0$, $\one_{L_1}$ euler, and with
precisely  one of the numbers in (\ref{e:aknum}) nonzero.
For instance  $(A) \vee (B)$ has only the first
 characteristic number $\sum_{p\in L_0} \beta (p)\ahlink \beta_2 (p)$
 nonzero.  We leave to the reader the construction of the remaining three
 examples.

\subsubsection {Elementary blocks, the general case}
We will define below an elementary block corresponding to
each function of the set $S=S_1\cup S_2\cup S_3\cup S_4$.
For reasons of
 symmetry, instead of $S_1$, $S_3$, we consider the
families
\begin{gather*}
 \alpha, \beta,\, \gamma; \tag{$S_1'$}\\
 \alpha_2, \, \beta_2 ,\,  \gamma_2, \alpha_3, \, \beta_3, \,\gamma_3 ;
\tag{$S_3'$}
\end{gather*}

We will use the elementary blocks to
construct $(L_1, \alpha, \beta, \gamma)$ and a vertex $p\in L_1$ with
given values of the following 31 expressions:
\begin{gather*}
\alpha (p), \beta (p),\, \gamma (p); \tag{$V_1$}\\
\ahlink (\alpha \beta)  (p),\, \ahlink (\alpha\gamma)  (p),\,
\ahlink  (\beta \gamma) (p),
\ahlink (\alpha \beta \gamma) (p);
\tag{$V_2$}\\
\ahlink \alpha_2  (p),\,  \ahlink  \beta_2
(p),\, \ahlink \gamma_2
 (p),\, \ahlink \alpha_3  (p),\, \ahlink  \beta_3 (p), \,\ahlink \gamma_3
(p);
\tag{$V_3$}
\end{gather*}
\begin{gather*}\tag{$V_4$}
\ahlink  (\alpha \beta_2) (p) ,\, \ahlink(\alpha \gamma_2) (p), \,
 \ahlink (\beta\alpha_2)  (p),\, \ahlink(\beta \gamma_2) (p), \,
\ahlink  (\gamma\alpha_2) (p) ,\, \ahlink(\gamma \beta_2) (p), \\
  \ahlink (\alpha \beta_3) (p),\,\ahlink ( \alpha \gamma_3) (p),  \,
 \ahlink (\beta\alpha_3) (p),\, \ahlink (\beta\gamma_3) (p),  \,
 \ahlink (\gamma\alpha_3) (p),\, \ahlink (\gamma\beta_3) (p),  \\
 \ahlink (\alpha\beta \gamma_2) (p),\,
\ahlink (\alpha \gamma \beta_2) (p),\,
\ahlink (\beta \gamma \alpha_2) (p), \\
 \ahlink(\alpha\beta \gamma_3) (p), \,
\ahlink (\alpha \gamma \beta_3) (p),\,
\ahlink (\beta \gamma \alpha_3) (p).
\end{gather*}

We will use the linear ordering of the set
$V=V_1\cup V_2\cup V_3\cup V_4$ given by the above list.
For each expression $v$ in $V$ we will define an elementary block
such that $v$ is nonzero at all vertices of the block, and all  expressions
which follow $v$ in the given ordering of $V$ vanish at all vertices of
the block. (The reader may easily check this property of the ordering of
$V$.)

The elementary blocks are defined by simple pictures.
On each picture we indicate the values of $(\alpha, \beta,\gamma)$ on
the 1-dimensional strata.  We illustrate 9 cases; the
remaining 22 elementary blocks
are defined by symmetry among the functions $\alpha$, $\beta,$ and $\gamma$.

\unitlength = 1cm
\begin{picture}(4,3)
\thicklines
\put(3,1){\circle{4}}
\put(2.3,1){\circle*{.2}}
\put(3.7,1){\circle*{.2}}
\put(0,0.9) {\shortstack{$(\alpha)$}}
\put(2.4,1.9) {\shortstack{$(1,0,0)$}}
\put(2.4,-.1) {\shortstack{$(1,0,0)$}}
\put(3.9,1) {\shortstack {$p$}}
\put(7,1) {\shortstack{$\alpha (p) \ne 0$}}
\end{picture}

\begin{picture}(4,3)
\thicklines
\put(3,1){\circle{4}}
\put(2.3,1){\circle*{.2}}
\put(3.7,1){\circle*{.2}}
\put(0,0.9) {\shortstack{$(\alpha_2)$}}
\put(2.4,1.9) {\shortstack{$(1,0,0)$}}
\put(2.3,-.1) {\shortstack{$(-1,0,0)$}}
\put(3.9,1) {\shortstack{$p$}}
\put(7,1) {\shortstack{$\ahlink \alpha_2 (p) \ne 0$}}
\end{picture}

\begin{picture}(4,3)
\thicklines
\put(3,1){\circle{4}}
\put(2.3,1){\circle*{.2}}
\put(3.7,1){\circle*{.2}}
\put(0,0.9) {\shortstack{$(\alpha_3)$}}
\put(2.4,1.9) {\shortstack{$(2,0,0)$}}
\put(2.4,-.1) {\shortstack{$(0,0,0)$}}
\put(3.9,1) {\shortstack{$p$}}
\put(7,1) {\shortstack{$\ahlink \alpha_3 (p) \ne 0$}}
\end{picture}



\unitlength = 1cm
\begin{picture}(4,3)
\thicklines
\put(4,1){\circle*{.2}}
\qbezier (4,1),(4.8,1.7),(5.6,1)
\qbezier (4,1),(4.8,.3),(5.6,1)
\qbezier (4,1),(3.9,2.3),(3.1,2.3)
\qbezier (4,1),(3.15,.9),(3.1,2.3)
\qbezier (4,1),(3.9,-.3),(3.1,-.3)
\qbezier (4,1),(3.15,1.1),(3.1,-.3)
\put(0,0.9) {\shortstack{$(\alpha  \beta)$}}
\put(1.8,1.7) {\shortstack{$(1,0,0)$}}
\put(1.8,-.1) {\shortstack{$(0,1,0)$}}
\put(4.3,.3) {\shortstack{$(1,1,0)$}}
\put(4.1,1.4) {\shortstack{$p$}}
\put(7,1) {\shortstack{$\ahlink (\alpha \beta) (p) \ne 0$} }
\end{picture}

\begin{picture}(4,4)
\thicklines
\put(4,1){\circle*{.2}}
\qbezier (4,1),(3.9,2.3),(2.8,2.3)
\qbezier (4,1),(3.15,.9),(2.8,2.3)
\qbezier (4,1),(3.9,-.3),(2.8,-.3)
\qbezier (4,1),(3.15,1.1),(2.8,-.3)
\qbezier (4,1),(4.1,2.3),(5.2,2.3)
\qbezier (4,1),(4.85,.9),(5.2,2.3)
\qbezier (4,1),(4.1,-.3),(5.2,-.3)
\qbezier (4,1),(4.85,1.1),(5.2,-.3)
\put(0,0.9) {\shortstack{$(\alpha  \beta  \gamma)$}}
\put(1.6,1.7) {\shortstack{$(1,0,0)$}}
\put(1.6,-.1) {\shortstack{$(0,1,0)$}}
\put(5.2,1.7) {\shortstack{$(0,0,1)$}}
\put(5.3,-.1) {\shortstack{$(1,1,1)$}}
\put(4.2,1.2) {\shortstack{$p$}}
\put(7,1) {\shortstack{$\ahlink (\alpha \beta \gamma) (p) \ne 0$}}
\end{picture}

\begin{picture}(4,3)
\thicklines
\put(3,1){\circle{4}}
\put(2.3,1){\circle*{.2}}
\put(3.7,1){\circle*{.2}}
\put(0,0.9) {\shortstack{$(\alpha \beta_2)$}}
\put(2.4,1.9) {\shortstack{$(1,1,0)$}}
\put(2.3,-.1) {\shortstack{$(1,-1,0)$}}
\put(3.9,1) {\shortstack {$p$}}
\put(7,1) {\shortstack{$\ahlink (\alpha \beta_2) (p) \ne 0$}}
\end{picture}

\begin{picture}(4,3)
\thicklines
\put(3,1){\circle{4}}
\put(2.3,1){\circle*{.2}}
\put(3.7,1){\circle*{.2}}
\put(0,0.9) {\shortstack{$(\alpha \beta_3)$}}
\put(2.4,1.9) {\shortstack{$(1,2,0)$}}
\put(2.4,-.1) {\shortstack{$(1,0,0)$}}
\put(3.9,1) {\shortstack{$p$}}
\put(7,1) {\shortstack{$\ahlink (\alpha \beta_3) (p) \ne 0$}}
\end{picture}

\begin{picture}(4,3)
\thicklines
\put(3,1){\circle{4}}
\put(2.3,1){\circle*{.2}}
\put(3.7,1){\circle*{.2}}
\put(0,0.9) {\shortstack{$(\alpha \beta \gamma_2)$}}
\put(2.4,1.9) {\shortstack{$(1,1,1)$}}
\put(2.3,-.1) {\shortstack{$(1,1,-1)$}}
\put(3.9,1) {\shortstack {$p$}}
\put(7,1) {\shortstack{$\ahlink (\alpha \beta \gamma_2) (p) \ne 0$}}
\end{picture}

\begin{picture}(4,3)
\thicklines
\put(3,1){\circle{4}}
\put(2.3,1){\circle*{.2}}
\put(3.7,1){\circle*{.2}}
\put(0,0.9) {\shortstack{$(\alpha \beta \gamma_3)$}}
\put(2.4,1.9) {\shortstack{$(1,1,2)$}}
\put(2.4,-.1) {\shortstack{$(1,1,0)$}}
\put(3.9,1) {\shortstack{$p$}}
\put(7,1) {\shortstack{$\ahlink (\alpha \beta \gamma ) (p) \ne 0$}}
\end{picture}

\subsubsection {Additivity, the general case.}
We want the expressions in $V$ to be additive with respect
to taking the wedge of elementary blocks.  This always
holds for 1-dimensional stratified sets and the expressions
$V_1$, $V_3$, $V_4$, but it can fail for $V_2$.  For the expressions in
$V_1$, $V_3$, $V_4$, additivity can be checked directly.
(For $\beta(p)$ additivity follows from (\ref{e:betas}), and
similarly $\alpha(p)$ and $\gamma(p)$ are additive. The expressions
in $V_3\cup V_4$ are of the form $\ahlink(\psi)(p)$, where $\psi$
is even, so $\ahlink(\psi)(p)\equiv \hlink(\psi)(p)\pmod 2$, and the latter
is clearly additive.) However, consider the expression
$\ahlink (\alpha \beta)$ in $V_2$, for instance.
The wedge $(L,p) = (L',p')\vee (L'', p'')$
of two 1-dimensional stratified sets
$(L', \alpha',\beta', \gamma')$
and $(L'',\alpha'', \beta'', \gamma'' )$ gives
\begin{equation*} 
\begin{split}
 \ahlink (\alpha  \beta) (p)  = & (\alpha  \beta) (p) - \hlink (\alpha
\beta) (p) \\
= & \alpha'  \beta' (p') + \alpha''  \beta'' (p'') +
(\alpha'(p') \cdot \beta'' (p'')
+ \alpha''(p'') \cdot \beta' (p')) \\
& -
(\hlink (\alpha' \beta') (p') + \hlink (\alpha'' \beta'') (p'') ) \\
= & \ahlink (\alpha' \beta') (p') + \ahlink (\alpha'' \beta'') (p'') +
(\alpha'(p') \cdot \beta'' (p'')
+ \alpha''(p'') \cdot \beta' (p')) ,
\end{split}
\end{equation*}
since $\hlink (\alpha  \beta)$ is clearly additive. There are similar
formulas for $\ahlink(\alpha\gamma)$, $\ahlink(\beta\gamma)$, and
$\ahlink(\alpha\beta\gamma)$.

Thus, even though the expressions in $V_2$ are not additive with
respect to taking the wedge, we still have that additivity holds
for them if one of the wedge factors is
 one of the elementary blocks $(\alpha  \beta) ,\, (\beta \gamma) ,\,
(\alpha \gamma ), (\alpha \beta \gamma)$,
since $\alpha =\beta=\gamma =0$ at the vertex of each of these blocks.

\subsubsection{Wedging elementary blocks.}

\begin{prop}\label{p:cons}
Given a subset $U$ of $V$,  there exists a compact
1-dimensional stratified set $(L_1,\alpha, \beta,\gamma)$ with a
distinguished vertex $p_0\in L_1$ such that
\begin{enumerate}
\item
    $\ahlink \alpha= \ahlink \beta =\ahlink \gamma = 0$, and
    $\alpha \beta$, $\alpha
        \gamma$, $\beta \gamma$, $ \alpha \beta \gamma$, and 
    $\one_{L_1}$ are euler;
\item
  $v(p_0) \ne 0$ if and only if $v\in U$;
\item
At every vertex $p\in L_1$, $p\ne p_0$, only one of the expressions
in $V$ is non-zero, and it is not an expression in  $V_2$.
\end{enumerate}
\end{prop}

\begin{rem} Note that the sum over all $p\in L_1$ of an expression
in $V_1$, $V_3$, or $V_4$ must be even.
(Since all functions of $S_3$ and $S_4$ are even,
 for each such a function $\psi$ we have that $\sum_{p\in L_0}
\ahlink \psi (p) = \int\ahlink \psi \dchi = \int \psi \dchi$ is
even.  Similarly the sums $\sum_{p\in L_0}
\alpha(p)$, $\sum_{p\in L_0} \beta(p)$, $\sum_{p\in L_0}
\gamma (p)$ are even.) This shows that the presence of nonzero
expressions at vertices other than $p_0$ is necessary
in the statement of the proposition.
\end{rem}

\begin{proof}[Proof of Proposition \ref{p:cons}]
We construct $L_1$ by wedging elementary blocks. By  definition of the
elementary blocks, condition (1) is
satisfied for any set constructed in this way.

The construction is by induction on the last element $u$ of $U$ with
respect to the given ordering of $V$.
   To fix ideas suppose that $u$ is in $V_3$ or $V_4$.
(The cases $u \in V_1$ or $u \in V_2$
are much simpler; they are left to the reader.)
Let $(u)$ be the elementary block corresponding to $u$. The underlying set
$B_u$ of this block has two vertices
$p_1$ and $p_2$, and $u$ is nonzero at both vertices.  It may happen that
$u$ is the only nonzero expression in $V$ at the vertices of $B_u$.
But in general
there are some expressions in $V$, necessarily prior to $u$ in the ordering
of $V$,
which are nonzero at both vertices of $B_u$.  Denote the set of such
expressions by $V_u$. By the inductive assumption applied to $V_u$,
there exists $(L',\alpha', \beta',\gamma')$ satisfying conditions (2) and
(3) for
$V_u$, with distinguished vertex $p'_0$.  Now wedge $(B_u,p_1)$
and $(L',p'_0)$ to obtain a set with all the
expressions of $V_u$ vanishing at the vertex $p_0=p_1\vee p'_0 $,
except perhaps for some expressions in
$V_2$. But any expression $v\in V_2$ which does not vanish at this
vertex can be corrected by
wedging with the correxponding elementary
block $(v)$.  Thus we can ensure that
all expressions in $V$ except $u$ vanish
at $p_0$.  The same procedure can be applied to the vertex
$p_2$ of $(u)$. In this way we
obtain a set $(L_u,\alpha_u,\beta_u,\gamma_u)$ with distinguished
vertex $p_0$ satisfying conditions (2) and (3) for the subset
$\{u\}$ of $V$.

Now let $U'=U\setminus\{u\}$. If $U'$ is empty we are done.
Otherwise we apply
the inductive assumption  to $U'$ and wedge the set
$(L',\alpha',\beta',\gamma')$
 we obtain with the set
$(L_u,\alpha_u,\beta_u,\gamma_u)$.  In this
way we can ensure that the  nonzero expressions at the distinguished
vertex $p_0$ are exactly those
of $\{u\}\cup U' = U$.  More precisely, we must again correct the
expressions $v\in V_2$ which do not vanish at $p_0$, by wedging with the
corresponding blocks $(v)$.
\end{proof}

\begin{cor}
Given a characteristic number $\cn (\mathbf m,\mathbf n)$,
there exists a compact stratified set $L_1$ of dimension 1 and
constructible functions $\alpha$, $\beta$,
$\gamma$ on $L_1$ such that
        \begin {enumerate}
        \item
        $\ahlink \alpha = \ahlink \beta =\ahlink \gamma = 0$ and $\alpha
\beta$,
$\alpha
        \gamma$, $\beta \gamma$, $\alpha \beta \gamma$, and
        $\one_{L_1}$ are euler;
        \item
        $\sum_{p\in L_0} \Phi(\mathbf m,\mathbf n)(p)$
        is the only nonzero sum of type (\ref {e:inv}).
        \end{enumerate}
\end{cor}
\begin{proof}
Recall that $\Phi(\mathbf m,\mathbf n) = \alpha^{m_\alpha}
\beta^{m_\beta}\gamma^{m_\gamma}\prod_{\psi\in S'}(\ahlink(\psi))^{n_\psi}$,
where the exponents are 0 or 1 and $\sum_{\psi\in S'}n_\psi>0$
(\ref{e:prod}).
Let $U$ be the subset of $V$ consisting of those expressions which
correspond to nontrivial factors of $\Phi(\mathbf m, \mathbf n)$. If we
apply Proposition
\ref{p:cons} to $U$, we obtain $(L_1,\alpha,\beta,\gamma)$ satisfying (1)
and such that $\sum_{p\in L_0} \Phi(\mathbf m, \mathbf n)(p)
\not\equiv 0 \pmod 2$. But $(L_1,\alpha,\beta,\gamma)$ has other nonzero
characteristic
numbers, namely those given by divisors of $\Phi(\mathbf m,\mathbf n)$.
 More precisely, if
$\mathbf m'=(m'_\alpha,m'_\beta,m'_\gamma)$ and $\mathbf n' =
(n'_\psi)_{\psi\in S'}$, with $m'_\alpha$, $m'_\beta$, $m'_\gamma$,
$n'_\psi$ equal to 0 or 1, $\sum_{\psi\in S'}n'_\psi>0$, and
$m'_\alpha\leq m_\alpha$, $m'_\beta\leq m_\beta$, $m'_\gamma \leq
m_\gamma$, $n'_\psi\leq n_\psi$ for all $\psi\in S'$, then
$\sum_{p\in L_0}\Phi(\mathbf m',\mathbf n')(p)\not\equiv 0\pmod 2$.
By induction on
the number of nontrivial factors of $\Phi(\mathbf m, \mathbf n)$ we
may assume that, for each such $(\mathbf m', \mathbf n')$ with
$(\mathbf m',\mathbf n')\neq(\mathbf m,\mathbf n)$, there exists
$(L_1',\alpha',\beta', \gamma')$ satisfying (1) and (2), so by taking
the disjoint union of $(L_1,\alpha,\beta,\gamma)$ with all such
$(L',\alpha',\beta',\gamma')$ we obtain the desired result.
\end{proof}

\subsection{Step 2. Construction of $(L_2, \varphi)$. }
\label{Step2}

Suppose we are given $(L_1, \alpha,\beta, \gamma)$ such that
$\ahlink \alpha = \ahlink \beta =\ahlink \gamma = 0$ and
$\alpha \beta, \alpha
\gamma, \beta \gamma$, $ \alpha \beta \gamma$, and $\one_{L_1}$ are euler.
  Note that the sets
$(L_1, \alpha,\beta, \gamma)$ constructed in step 1 (section \ref{Step1})
 are all are actually constructed by wedging blocks of the type

\begin{picture}(4,3,2)
\thicklines
\put(4,1.5){\circle{4}}
\put(3.3,1.5){\circle*{.2}}
\put(4.7,1.5){\circle*{.2}}
\put(0,1.4) {\shortstack{$(P,\alpha,\beta, \gamma) \quad = $}}
\put(3.2,2.4) {\shortstack{$(\alpha', \beta',\gamma')$}}
\put(3.1,.4) {\shortstack{$(\alpha'', \beta'',\gamma'')$}}
\end{picture}

\noindent
where $\alpha'\equiv \alpha'' \pmod 2$, $\beta' \equiv \beta''
\pmod 2$, and $\gamma'\equiv \gamma'' \pmod2$.

\subsubsection{Cycles and the values of $\alpha$, $\beta$, $\gamma$
mod 2.}
Let $(L_2, \varphi)$ be such that $\hlink
\varphi =0$, and let $S\subset L_2$ be a 1-dimensional stratum.
Suppose that $S$ is in the boundary
of exactly $k$
two-dimensional strata $S_1, \ldots, S_k$, with the
values of $\varphi$ on the strata equal to $l_1, \ldots,l_k$, respectively.
Let $p\in
S$.  Then, since $\hlink \varphi = 0$,
\begin{equation} \label {e:betas2}
\begin{split}
& \alpha(p) = \varphi (p) = \half \sum l_i, \\
& \beta (p) = \hlink \varphi^2(p) = \varphi ^2 (p)
- \half \sum l^2_i, \\
& \gamma (p) =  \hlink \varphi^3(p) = \varphi ^3 (p)
- \half \sum l^3_i ,
\end{split}
\end{equation}
 which are formulas analogous to \eqref{e:betas}.

Consider a block $(P, \alpha, \beta,
\gamma)$ as pictured above.  By the computations of
(\ref{blockexample}) there
exist $l_1,\dots,l_k$ such that  \eqref{e:betas2} gives $(\alpha, \beta,
\gamma)$ modulo $2$.
(The value $(\alpha,\beta,\gamma)\equiv(1,0,0)\pmod 2$ is given by
$(l_1,l_2)=(1,1)$;
$(\alpha,\beta,\gamma)\equiv(0,1,0)\pmod 2$ is given by $(l_1,l_2)=(1,-1)$;
$(\alpha,\beta,\gamma)\equiv(1,1,1)\pmod 2$ is given by $(l_1,l_2)=(2,0)$.)
Define $(L_2, \varphi)$ to be the disjoint
union of $k$ two-dimensional discs $D_i$, $i=1,\ldots, k$, identified
along their boundaries, and let $\varphi$ be equal to $l_i$ on the
interior of each disc $D_i$ and equal to $\half \sum l_i$ on their common
boundary.  Then  $(L_2, \varphi)$ satisfies some of the desired
properties for $(P, \alpha, \beta, \gamma)$; namely,
$\hlink \varphi = 0$ and
$\alpha \equiv \varphi |_{L_1}$, $\beta \equiv \hlink \varphi^2$,
$\gamma
\equiv \hlink \varphi^3 \pmod 2$ on open 1-dimensional strata.

We generalize this construction as follows.
\begin{prop}\label{p:cycles}
 Let $(L_1,\alpha,\beta,\gamma)$ be such that $\ahlink\alpha=
 \ahlink\beta=\ahlink\gamma=0$ and $\alpha\beta$, $\alpha\gamma$,
 $\beta\gamma$, $\alpha\beta\gamma$, and $\one_{L_1}$ are euler.
 There exists $(L_2,\varphi)$ such that
\begin{enumerate}
\item
$\hlink \varphi = 0$,
\item
  $\alpha \equiv \varphi |_{L_1}$, $\beta \equiv \hlink \varphi^2$,
$\gamma
\equiv \hlink \varphi^3 \pmod 2$ on open 1-dimensional strata.
\end{enumerate}
\end{prop}

\begin{proof}
It is easy to see the hypotheses imply that, for all
$(a,b,c)\in(\Z/2)^3$, the set $L(a,b,c)=\{p\in L_1\ |\
(\alpha,\beta,\gamma)\equiv (a,b,c)\pmod 2\}$ is a 1-cycle (mod 2).
We construct $(L_2,\varphi)$ by induction on the number of edges
(open 1-strata) of $L_1$ on which $(\alpha,\beta,\gamma)$ is
nontrivial (mod 2). Suppose $(\alpha,\beta,\gamma)\equiv(a,b,c)
\not\equiv(0,0,0)\pmod 2$ on some edge of $L_1$. Since $L(a,b,c)$ is
a cycle (mod 2), there is a subset $C(a,b,c)$ of $L_1$ such that
$C(a,b,c)$ is a union of closures of edges of $L_1$, $C(a,b,c)$ is
homeomorphic to a circle, and
$(\alpha,\beta,\gamma)\equiv(a,b,c)\pmod 2$ on every edge of $C(a,b,c)$.
Suppose $(L_2',\varphi')$ satisfies (1) and (2) for the set
$L_1'$ obtained by removing the edges of $C(a,b,c)$ from $L_1$.
Then we obtain $(L_2,\varphi)$ from $(L_2',\varphi')$ by attaching
discs along $C(a,b,c)$ just as we did for $(P,\alpha,\beta,\gamma)$
above.
\end{proof}

\subsubsection {Wedges, bubbles, and the values of $\alpha$, $\beta$,
$\gamma$ mod $4$.} \label {bubbles}

Suppose we are given two 2-dimensional stratified sets
$(L',\varphi')$, $\hlink
\varphi'=0$,
and $(L'',\varphi'')$, $\hlink \varphi''=0$, with distinguished
1-dimensional strata $S'\subset L'$ and $S''\subset L''$.
Identifying $S'$ and $S''$ we produce a new 2-dimensional set
$L$ with distinguished 1-stratum $S$.
(If the boundary of $S'$ in $L'$ is a single point $p$ and the
boundary of $S''$ in $L''$ is two points $p_1$ and $p_2$, then
both $p_1$ and $p_2$ are identified with $p$.)
The functions $\varphi'|_{L'\setminus \bar S'}$ and
$\varphi''|_{L''\setminus \bar S''}$  extend to a
constructible function $\varphi$ on $L$ satisfying $\hlink \varphi = 0$.
The function $\varphi$ is uniquely defined except at the vertices on the
boundary of $S$.  We call $(L,\varphi,S)$ the {\it wedge} of
$(L',\varphi',S')$ and $(L'',\varphi'',S'')$,  and we denote it by
$(L',\varphi',S') \vee (L'',\varphi'',S'')$.
The values of $\alpha = \varphi|_{L_1}$, $\beta =
\hlink \varphi^2$, $\gamma = \hlink \varphi^3$ on $S$, are
determined by the values of
$\alpha' = \varphi'|_{L'_1}$, $\beta' =
\hlink (\varphi')^2$, $\gamma' = \hlink (\varphi')^3$,  and
$\alpha'' = \varphi''|_{L''_1}$, $\beta'' =
\hlink (\varphi'')^2$, $\gamma'' = \hlink (\varphi'')^3$, on $S'$ and $S''$,
respectively, by
\begin{equation} \label {e:add2}
\begin{split}
& \alpha = \alpha' + \alpha'' , \\
& \beta = \beta' + \beta'' + 2 \alpha' \cdot \alpha'' , \\
& \gamma =  \gamma'+ \gamma'' +3\alpha' \cdot \alpha''
(\alpha' + \alpha'') .
\end{split}
\end{equation}

In particular, as we have noticed already for the
analogous situation in section \ref{Step1},
$(\alpha, \beta, \gamma) \pmod 2$ is additive with respect to the
operation of taking the wedge.  On the other hand, it is not in
general additive if we consider the values mod 4.

Fix $l\in \Z$.
Let $B_l$ denote the 2-dimensional sphere stratified by the eastern
and western hemispheres, the $0^\circ$ and $180^\circ$ meridians, and
the north
and south poles.  Set $\varphi$ on $B_l$ to be constant and equal
to $l$.   The $0^\circ$ meridian will be the distinguished
1-dimensional stratum $S_l$ of $B_l$. The triple $(B_l,\varphi, S_l)$
will be called a {\it bubble}.

\begin{prop}\label{p:bubbles}
Let $(L_2,\varphi)$, be such that $\hlink \varphi =0$, and let $S$ be a
1-dimensional stratum of $L_2$. Let $ \alpha_0$, $\beta_0$, $\gamma_0$
denote the values of $\varphi$, $\hlink \varphi^2$, $\hlink \varphi^3$
on $S$.  Let $\alpha'_0 \equiv \alpha_0 \pmod 2$,
$\beta'_0 \equiv \beta_0 \pmod 2$, $\gamma'_0 \equiv \gamma_0 \pmod 2$.
Then there exist $l_1,\dots,l_k$, such that the
wedge of $(L_2,\varphi,S)$ along $S$ with the bubbles
$(B_{l_i},\varphi,S_{l_i})$
 has the values
$(\alpha,\beta,\gamma)\equiv ( \alpha'_0, \beta'_0, \gamma'_0) \pmod 4$ on
the distinguished
stratum of the wedge.
\end{prop}

\begin{proof}
When we wedge the bubbles $B_l$ along their distinguished strata we may
get any combination of even values of $\alpha$, $\beta$, $\gamma$ modulo 4.
In particular we have:
\begin{equation*}
\begin{split}
 B_1\vee B_3\  \text { gives }  (\alpha,\beta, \gamma) = (0,2,0) ,
 \\
 (\vee_3 B_1)\vee B_2 \vee B_3\   \text { gives }
(\alpha,\beta, \gamma) = (0,0,2) ,   \\
 (\vee_4 B_1)\vee B_2\  \text { gives }   (\alpha,\beta, \gamma) = (2,0,0).
\end{split}
\end{equation*}
So the proposition follows from the formulas \eqref{e:add2}.
\end{proof}

\begin{cor} \label{c:step2}
Given $(L_1, \alpha,\beta, \gamma)$ such that
$\ahlink \alpha = \ahlink \beta =\ahlink \gamma = 0$ and
$\alpha \beta$, $\alpha \gamma$, $\beta \gamma$, $ \alpha \beta
\gamma$, and $\one_L$ are euler, then there exists
$(L_2,\varphi)$ satisfying the claim of step 2.
\end{cor}

\begin{proof}
Propositions \ref{p:cycles} and \ref{p:bubbles} guarantee
the existence of $(L_2,\varphi)$ such that
\begin {enumerate}
        \item
        $\hlink \varphi =0$.
        \item
        $\alpha \equiv \varphi |_{L_1}, \beta \equiv \hlink \varphi^2,
        \gamma
        \equiv \hlink \varphi^3 \pmod 4$ on open 1-dimensional strata.
        \end{enumerate}
We establish 2 (mod 2) at the vertices.  First note that we may
arbitrarily change
the value of $\varphi$ at the vertices without
destroying property 1.   Hence we may arbitrarily adjust the
values of $\varphi$ so that 2 holds for $\alpha$.  Then the
analogous property for $\beta$ and $\gamma$ is satisfied
automatically.  Indeed,
the difference $\psi = \beta - \hlink
\varphi^2$ is supported in $L_1$ and
is divisible by $4$ on open 1-strata.  Therefore the property
 $\ahlink \psi =0$ gives that $\psi$ is divisible by $2$ on
$L_0$, as required. A similar argument works for $\gamma$.

Finally if we wedge $L_2$ with a 2-sphere $S^2$ at vertices $p'\in
L_2$ and $p''\in S^2$, with $\varphi=1$ on $S^2$, then we change the
integral of $\varphi$
without changing any of the previous properties.
This completes the proof of the corollary and step 2.
\end{proof}

\subsection{Step 3. Construction of $L_3$. }
\label{Step3}

It suffices to prove the following result.

\begin{prop}\label{p:step3}
Let $\varphi$ be a constructible function defined on a compact
semialgebraic set $Y$ of dimension $\le 2$, and suppose $\Lambda
\varphi = 0$.  There exists a compact semialgebraic set
$\tilde Y$ of dimension $\le 3$, such that $Y\subset \tilde Y$,
$\tilde Y$ is euler, $\supp  \ahlink \one_{\tilde Y} \subset Y $,
and $\varphi \equiv \ahlink \one_{\tilde Y} \pmod {2\ideal (Y)}$.
\end{prop}

\begin{proof}
Let $T$ be a triangulation of $Y$ such that $\varphi$ is constant
on open simplices of $T$.  For each simplex $\Delta$ of $T$ we
construct a set $\tilde \Delta$, with $\Delta \subset \tilde
\Delta$.
 Then $\tilde Y$ will be defined as $Y$ with each $\tilde \Delta$
 attached along the simplex $\Delta$, for all $\Delta \in T$.
 In other words, if $Z$ is the disjoint union of all the sets
 $\tilde\Delta$ for $\Delta\in T$, then $\tilde Y$ is obtained
 from the disjoint union of $Y$ and $Z$ by identifying each point
 $y\in Y$ with its image in $\tilde\Delta$, for all $\Delta\in T$
 containing $y$.

Given
$\Delta \in T$, let $\varphi (\Delta)$ be the value of
$\varphi$ on the interior of $\Delta$.   The construction of
$\tilde \Delta$ depends on the dimension of  $\Delta$.

$\dim \Delta = 0$:
Then $\Delta$ is a one point space say $\Delta = \{p\}$.  Let
$m= m(\Delta)$
 be a positive integer such that $m \equiv 1 - \varphi(\Delta) \pmod
2$.  Then $\tilde \Delta$ is defined to be the wedge of $m$
circles at $p$.

$\dim \Delta = 1$:
Let $m = m(\Delta)$ be a positive integer such that
$m \equiv  \varphi(\Delta) \pmod 4$.
Then $\tilde \Delta$ is defined to be the union of $m$
segments with their left endowment identified and their right
endowment identified; that is, $\tilde \Delta$  is homeomorphic to
the suspension of $m$ points.  We identify $\Delta$ with one of
the segments.

$\dim \Delta = 2$:
Let $m= m(\Delta)$ be a positive integer such that
$m \equiv 1 - \varphi(\Delta)  \pmod 8$.
 Let $R_m$ denote the wedge of m circles and let
$x_0\in R_m$ denote the base point of the wedge.
Then $\tilde \Delta$ is defined to be $\Delta \times
R_m/\sim$, where we collapse $\{p\}\times R_m$ for each $p\in
\partial \Delta$;  that is, $(p,x)\sim (q,y)$ if
and only if $p=q\in \partial
\Delta$ or $p=q \in \operatorname{int}(\Delta)$ and $x=y$.  $\Delta$ is
identified with $\Delta\times \{x_0\}$.

Let  $\tilde Y$ be the union of $Y$ and the sets $\tilde \Delta$
 attached along the simplices $\Delta$.  Then, since $\tilde Y
 \setminus Y$ is a union of open 1-cells and open 3-cells,
 $\supp \ahlink \one_{\tilde Y} \subset Y$.
Let $\Delta^2\in T$ be a
simplex of dimension 2, and let $p\in \operatorname{int}(\Delta^2)$.  Then
\begin{equation} \label{e:dim2}
 \ahlink \one_{\tilde Y} (p) = 1 - m(\Delta^2) \equiv \varphi(p)
\pmod 8 ,
\end{equation}
as required.

Let $\Delta^1\in T$ be a
simplex of dimension 1, and let $p\in \operatorname{int}(\Delta^1)$.
Let $W$ be a
normal slice to  $\Delta^1$ in $\tilde Y$.
(In other words, $W=H\cap\tilde Y$, where $Y$ is realized as a
semialgebraic set in $\R^n$ and $H$ is an $(n-1)$-plane in $\R^n$
which meets $\Delta^1$ transversally at an interior point of
$\Delta^1$.)
  Let
$\psi = \varphi -  \ahlink \one_{\tilde Y}$.  Then $\Lambda \psi
= 0$, which gives, by the slice property (section \ref{properties}(d))
  $\Omega (\psi|_H)= 0$; that is,
\begin{equation*}
2\psi (p) - \Lambda (\psi|_H) (p) = 0 .
\end{equation*}
 By \eqref{e:dim2}, $\Lambda (\psi|_H) (p)\equiv 0 \pmod 8$,
and consequently $\psi (p) \equiv 0 \pmod 4$ as required.

Finally, let $p$ be a vertex of $X$.  Let $L$, resp. $\tilde L$,
 denote the link of
$Y$, resp. $\tilde Y$,  at $p$.  The intersections of $Y$ and
 $\tilde Y$ with a small  sphere centered at $p$ induce cell
 structures on  $L$ and $\tilde L$.
Denote by $C^0$ and $C^1$, the cells of $L$ of dimension $0$
and $1$, respectively.
Let $\hat \varphi = \varphi|_L$.   For a cell $C$ we denote by
 $\hat \varphi (C)$ the value of
$\hat \varphi$ on the interior of $C$.  By the
slice property (section \ref{properties}(d)),
$\Omega \hat \varphi = 0$; that is, for each
0-cell $C^0$ in $L$,
$\hat \varphi (C^0) = \half \sum_{C^1 \supset C^0}
\hat \varphi (C^1)$.
Consequently
\begin{equation} \label {e:2hats}
\sum_{C^0\subset L} \hat \varphi (C^0) = \sum_{C^1\subset L}
\hat \varphi (C^1) .
\end{equation}

By construction,
\begin{equation}\label {e:3hats}
\begin{split}
& \# \{\text { 0-cells in $\tilde L$ }\} \equiv 2(1-\varphi (p))
+ \sum_{C^0\subset L} \hat \varphi (C^0) \pmod 4 \\
& \# \{\text { 1-cells in $\tilde L$ }\} \equiv
\# \{\text { 1-cells in $L$ }\}  \pmod 4 \\
& \# \{\text { 2-cells in $\tilde L$ }\} \equiv
\sum_{C^1\subset L} (1-\hat \varphi (C^1)) \pmod 4
\end{split}
\end{equation}
By \eqref{e:2hats} and \eqref{e:3hats},
$\chi (\tilde L) \equiv  2(1-\varphi (p)) \pmod 4$, so
\begin{equation*}
\ahlink \one_{\tilde Y} (p) = 1 - \half \chi (\tilde L) \equiv
\varphi(p) \pmod 2 ,
\end{equation*}
 as required.
\end{proof}

\medskip

\section{Arithmetic operators on algebraically constructible
  functions} \label{Arith}

\subsection{Polynomial operators with rational coefficients.}
\label{Pclass}

Since the set of algebraically constructible functions
$\algc (X)$ is a ring, each
polynomial $P$ with integer coefficients gives an operator
$\varphi
\to P(\varphi)$ on $\algc (X)$.  We note that some other
polynomials with rational coefficients have the same property.

\begin{lem}\label{l:4-2}
If $\varphi$ is an algebraically constructible function on $X$, then
$\half (\varphi^4-\varphi^2)$ is also algebraically
constructible.
\end{lem}

\begin{proof}
By \eqref{e:signs} we may suppose that
$\varphi = \sum_{i=1}^s \sign g_i$, where $g_i$
are polynomials.  Then
\begin{equation*}
\begin{split}
\varphi ^2 & =  \sum (\sign g_i)^2 + 2 \sum_{i<j} \sign (g_ig_j),
\\
\varphi ^4 & =  \sum (\sign g_i)^4 + 2 \sum_{i<j} \sign
(g_ig_j)^2 \\
&\ \ \ \ + 4  \big(\sum (\sign g_i)^2\big) \big( \sum_{i<j} \sign
(g_ig_j) \big) +
4  \big(\sum_{i<j} \sign (g_ig_j) \big)^2 ,
\end{split}
\end{equation*}
and hence the result.
\end{proof}

The existence of such operators on algebraically constructible
functions gives even more local conditions on the topology
  of real algebraic sets.  We will discuss these conditions for the
sets of dimension $4$ in the following subsections.  But first we
classify polynomials which give operations on $\algc(X)$ and show,
 in particular, that modulo $8$ the ring of
all such polynomials is generated by $\Z [t]$ and
  $P(t) = \half (t^4 - t^2)$.
First we note that every such polynomial
takes integer values on integers and hence can be uniquely
written as a sum
\begin{equation*}
P = \sum_{p\ge 0} n_p f_p ,
\end{equation*}
where $n_p\in \Z$, and $f_p(t)$ are the binomial
polynomials:
\begin{equation}\label{e:binom}
f_p (t) = \frac {t(t-1) \cdots (t-p+1)} {p!}, \quad p= 1,2, \ldots,
\quad f_0(t) = 1 .
\end{equation}
\medskip

\begin{thm}\label{t:pol}
Let $P\in \Q [t]$.  The following conditions are equivalent:
\begin{itemize}
\item[(i)]
For every real algebraic set $X$ and every algebraically constructible
function $\varphi$ on
$X$, $P(\varphi)$ is algebraically constructible;
\item[ (ii)]
$P$ preserves finite formal sums of signs; that is,
 for all $s$, all the coefficients $a_\alpha$,
 $\alpha=(\alpha_1,\ldots,\alpha_s)\in \{0,1,2\}^s$,  in  the
formal expansion
\begin{equation*}
P\big(\sum_{i=1}^s \sign g_i\big) = \sum_\alpha a_\alpha \prod_{i=1}^s
(\sign g_i)^{\alpha_i}
\end{equation*}
are integral;
\item[ (iii)]
$P = \sum_{p\ge 0} n_p f_p$ and each $n_p$ is an integer
divisible by $2^{[p/2]}$.
\end{itemize}
\end{thm}

\begin{proof}
Note that (ii) implies (i) easily by \cite{parszaf1},
\cite{parszaf2}, {\it i.e.}~by the characterization  \eqref{e:signs}
of algebraically constructible functions. On the other hand,
(i)$\Longrightarrow$(ii) is not so obvious.  We will
proceed in a different way and establish
(ii)$\Longleftrightarrow$(iii) first.

Denote by $\pol$ the set of polynomials satisfying (ii).
 It is easy to see that $\pol$ is a ring
invariant by translations,
{\it i.e.}~$P(t) \in \pol$ if and only if $P(t+1) \in \mathcal
P$.  Define
\begin{eqnarray*}
\Delta P (t) = P(t+1) - P(t)  .
\end{eqnarray*}

\begin{lem}\label{l:delta}
$P\in \pol$ if and only if  $P$ is integer-valued on integers,
$\Delta P\in \pol$, and $\half \Delta ( \Delta P) \in \pol$.
\end{lem}

\begin{proof}
Define
\begin{equation}
\Delta_2 P (t) = P(t+1) - P(t-1) .
\end{equation}
Then we may replace the condition $\half \Delta ( \Delta P) \in \pol$
in the statement of the lemma
 by
$\half \Delta_2 P \in \pol$.
Indeed,
\begin{equation}
\half \Delta_2 P (t+1) = \half (P(t+2) - P(t)) = \half
\Delta ( \Delta P)(t) +\Delta P (t) ,
\end{equation}
and hence $\Delta P\in \pol$ and $\half \Delta ( \Delta P) \in \pol$
if and only if $\Delta P\in \pol$ and $\half \Delta_2 P \in \pol$.

Write $P = \sum_{p\ge 0} n_p f_p$, $n_p\in \Z$, and consider the
coefficients $a_\alpha$,
$\alpha=(\alpha_0,\alpha_1,\ldots,\alpha_s)\in \{0,1,2\}^{s+1}$,
given by
\begin{equation*}
P\big(\sum_{i=0}^s \sign g_i\big) =
\sum_\alpha a_\alpha \prod_{i=0}^s (\sign g_i)^{\alpha_i}.
\end{equation*}
Note that $a_\alpha$ are not necessarily integers if $P\notin \pol$.
Setting $g_0 = 1,-1, \text { or } 0$, we obtain
\begin{gather*}
P\big(\pm 1 + \sum_{i=1}^s \sign g_i\big) =
\sum_{\alpha'}
(a_{0,\alpha'} \pm a_{1,\alpha'} + a_{2,\alpha'})
\prod_{i=1}^s (\sign g_i)^{\alpha_i}, \\
P\big(\sum_{i=1}^s \sign g_i\big) =
\sum_{\alpha'}
a_{0,\alpha'} \prod_{i=1}^s (\sign g_i)^{\alpha_i},
\end{gather*}
where $\alpha'= (\alpha_1,\ldots,\alpha_s)$.
Consequently,
\begin{equation*}
\Delta P\big(\sum_{i=1}^s \sign g_i\big) =
P\big(1 +\sum_{i=1}^s \sign g_i\big) - P\big(\sum_{i=1}^s \sign g_i\big)
 = \sum b_{\alpha'}\prod_{i=1}^s (\sign g_i)^{\alpha_i} ,
\end{equation*}
where
\begin{equation}\label{e:coef1}
 b_{\alpha'} =  a_{1, \alpha'} +  a_{2, \alpha'}.
\end{equation}
Similarly,
\begin{equation*}
\half \Delta_2 P\big(\sum_{i=1}^s \sign g_i\big) =
\sum_{\alpha'} c_{\alpha'}\prod_{i=1}^s (\sign g_i)^{\alpha_i} ,
\end{equation*}
where
\begin{equation}\label{e:coef2}
 c_{\alpha'} =  a_{1, \alpha'}.
\end{equation}
In particular, if all $a_\alpha$ are integers, so are all
$b_{\alpha'}$ and $c_{\alpha'}$.  This shows that if
$P\in \pol$ then $P$ is integer-valued on integers,
$\Delta P\in \pol$, and $\half \Delta ( \Delta P) \in \pol$.

Now we show the converse.  Suppose that all $b_{\alpha'}$ and
$c_{\alpha'}$ are integers.  Then all the coefficients $a_{1, \alpha'}$
and $a_{2, \alpha'}$ are also integral.  Since $a_\alpha$ are
invariant by permutations of $\alpha_0,\alpha_1,\ldots,\alpha_s$, this
implies that all $a_\alpha$, except maybe $a_{(0,\ldots,0)}$, are
integral.  But $a_{(0,\ldots,0)}= n_0= P(0)$ is integral by assumption.
This completes the proof of the lemma.
\end{proof}

To show (ii)$\Longleftrightarrow$(iii) we note that
$\Delta f_p = f_{p-1}$.
Clearly $f_0$ and $f_1$ are in $\pol$.  Since $\deg \Delta P < \deg
P, $ (ii)$\Longleftrightarrow$(iii)
 follows easily from Lemma \ref{l:delta} by induction
on the degree of $P$.

To complete the proof of the theorem we show (i)$\Longrightarrow$(iii).
Suppose $P = \sum_{p\ge 0} n_p f_p$ and suppose that there exists at
least one $p$ such that $2^{[p/2]}\nmid n_p$.  Let $p_0$ be the
smallest $p$ with this property, and let $k=[p_0/2]$.  Consider
the algebraically constructible function
$\varphi_{p_0} =\sign x_1 + \ldots +\sign x_k  + p_0 - k$ on
$\R^k$, where $x_1,\ldots,x_k$ denote the coordinate functions on $\R^k$.
Then $\varphi_{p_0}$ takes values between $0$ and $p_0$ for
$p_0$ even, and between $1$ and $p_0$ for $p_0$ odd.  Consider the function
\begin{equation}
P(\varphi_{p_0}) = \sum_{p<p_0} n_p f_p(\varphi_{p_0}) +
n_{p_0}f_{p_0}(\varphi_{p_0}) + \sum_{p>p_0} n_p f_p(\varphi_{p_0}) .
\end{equation}
The first summand is an algebraically constructible function, and the
last summand vanishes identically.  Note that
$f_{p_0}(\varphi_{p_0})$ is the characteristic function of the
open first quadrant $Q_k$ of $\R^k$.  So
(i)$\Longrightarrow$(iii) is a consequence of the following lemma.

\begin{lem}\label{l:2tok}
Let $Q_k$ be the open or closed first quadrant in $\R^k$.
For $n\in \Z$, the function $n\cdot\one_{Q_k}$ is
algebraically constructible if and only if $2^k | n$.
\end{lem}
\begin{proof}
Let $Q_k$ denote the closed first quadrant
and let
$Q'_k$ denote the open first quadrant
in $\R^k$. Suppose $\varphi= n\cdot
\one_{Q_k}$ is algebraically constructible. Let $\one_{\R^{k-1}}$ be
the characteristic function of $\R^{k-1}$, and let
$\psi = (\hlink\varphi)\one_{\R^{k-1}}$.
Then $\psi$ is algebraically constructible by Theorem \ref{t:hlink},
and $\psi = \frac n2 \one_{Q_{k-1}}$, with $\frac n2\in\Z$. By
induction on $k$ we have $2^{k-1}|\frac n2$, so $2^k|n$. If
$\varphi'= n\cdot\one_{Q'_k}$, and $\psi'=(\hlink\varphi')\one_{\R^{k-1}}$,
then $\psi'=(-1)^{k-1}\frac n2
\one_{Q_{k-1}}$, and so again $2^k|n$ by induction.

To prove the converse, let $S\subset\{1,\dots,k\}$, and let
$f_S(x_1,\dots,x_k)=(\epsilon_1x_1^2,\dots,\epsilon_kx_k^2)$, where
$\epsilon_i=1$ for $i\in S$ and $\epsilon_i=0$ for $i\not\in S$. Then
\begin{equation*}
\begin{split}
& 2^k\one_{Q_k}(x) = \sum_S\chi(f_S^{-1}(x)),\\
&2^k\one_{Q_k'}(x) = \sum_S(-1)^{k-|S|}\chi(f_S^{-1}(x)),
\end{split}
\end{equation*}
so these functions are algebraically constructible by \eqref{e:algcons}.
\end{proof}
 This completes the proof of Theorem \ref{t:pol}.\end{proof}

Let $\pol$ denote the ring of polynomials in $\Q [t]$ which
satisfy one of the equivalent conditions of Theorem \ref {t:pol}.
  Let $\pola$ denote the ring of polynomials in $\Q [t]$ which take
  integer values  on integers.
By condition (iii) of Theorem \ref{t:pol}, $\pol$ is additively
generated by $\Z [t]$ and
the polynomials in $\pol$ whose values are divisible by $4$;
that is, $\pol/(\pol \cap 4\pola)$ is additively generated by
$\Z[t]$.
This is the reason that polynomial operators on algebraically
constructible functions do not
give new characteristic numbers for real algebraic set of dimension $\le 3$.

On the other hand, Theorem \ref{t:pol} shows that $\pol/(\pol \cap
8\pola)$ is additively generated by $2^{[p/2]}f_p(t)$,
$p=0,\dots,5$.
Later on we will use the following equivalent system of additive
generators:
\begin{equation} \label{e:mod8}
\begin{split}
& 1,\ t,\ t^2 -t,\ t^3-t,\  P_4(t)= \half t(t-1)(t-2)(t-3), \\
&  P_5(t)= \half t(t-1)(t-2)(t-3)(t-4) .
\end{split}
\end{equation}
Note that, in particular, the polynomial of Lemma \ref{l:4-2} can be
written
\begin{equation}
\half (t^4-t^2) = \half t(t-1)(t-2)(t-3) + 3t(t-1)^2 .
\end{equation}

Theorem \ref{t:pol} admits an obvious generalization to the case
of polynomials of many variables.

\begin{cor}
Let $P\in \Q[t_1,\ldots t_s]$.  The following conditions are equivalent:
\begin{itemize}
\item[(i)]
For every real algebraic set $X$ and all algebraically constructible
functions $\varphi_1,\ldots, \varphi_s$ on
$X$, $P(\varphi_1, \ldots, \varphi_s)$ is algebraically constructible;
\item[ (ii)]
$P$ preserves finite formal sums of signs;
\item[ (iii)]
$P = \sum_{p=(p_1,\ldots, p_s)} n_p \prod_{i=1}^s  f_{p_i} (t_i)$
and each $n_p$ is an integer
divisible by $2^{[p_1/2]+\ldots +[p_s/2]}$.
\end{itemize}
\end{cor}

\begin{proof}
(ii)$\Longrightarrow$(i) follows again from \eqref{e:signs}. 
(iii)$\Longrightarrow$(ii) results directly from the analogous implication
in
Theorem \ref{t:pol}. The proof of (i)$\Longrightarrow$(iii) is similar
to the
proof of the corresponding implication in
Theorem \ref{t:pol}. We sketch the proof.   Suppose
$P = \sum_{p=(p_1,\ldots, p_s)} n_p \prod_{i=1}^s  f_{p_i} (t_i)$
has integral coefficients, and suppose that there exists at least one
$p$ such that $2^{[p_1/2]+\ldots +[p_s/2]}\nmid n_p$.
 Let $p_0$ be the
smallest $p$ (with respect to the total degree or lexicographic order)
 with this property, and let $k_i=[p_{0i}/2]$.  Consider
$P(\varphi_{p_{01}},\ldots ,\varphi_{p_{0s}} )$,
where as before
$\varphi_{p_{0i}} =\sign x_1 + \ldots +\sign x_{k_i}  + p_{0i} - k_i$
are algebraically constructible functions on
$\R^{k_i}$.  Then
\begin{equation}\label{e:P}
P(\varphi_{p_{01}},\ldots ,\varphi_{p_{0s}} ) =
\psi +
n_{p_0}\prod_{i=1}^s f_{p_{0i}}(\varphi_{p_{0i}}) ,
\end{equation}
where $\psi$ is algebraically constructible. Now
$\prod_{i=1}^s f_{p_{0i}}(\varphi_{p_{0i}})$ is
 the characteristic function of the
open first quadrant in $\R^{k_1}\times \cdots \times \R^{k_s}$
and hence, by Lemma \ref {l:2tok},  the function \eqref{e:P} is
not algebraically constructible.
\end{proof}

\subsection{Construction of invariants.}\label{Pcons}
Using Theorem
\ref{t:pol} we now define new depth $k$ euler conditions and
characteristic numbers.  Since our procedure goes along
the lines of section \ref{Cons}, we just outline the construction and
emphasize the new features.

\subsubsection{$\pol$-rings and $\pol$-euler spaces}\label
{Prings}
Let $\pol \subset \Q [t]$ be the set of polynomials described in
Theorem \ref{t:pol}.  Given a constructible function $\varphi$, we let
$\pol (\varphi) = \{P(\varphi)\, |\, P\in \pol\}$.  More generally,
if $\Phi$ is a set of constructible functions, we denote by
$\pol(\Phi)$ the ring generated by $\pol(\varphi)$ for $\varphi \in
\Phi$, {\it i.e.} the smallest ring containing $\Phi$ and closed under
the operators given by $P\in \pol$.  We call $\pol(\Phi)$ the
{\it $\pol$-ring} generated by $\Phi$.

Suppose now that the semialgebraic set $X$ is  homeomorphic to
an algebraic set.  Then by Theorems \ref{t:hlink} and \ref{t:pol} all
functions constructed from $\one_X$ by means of the arithmetic
operations $+$, $-$, $*$, and the operators $\hlink$ and $\pol$,
are integer-valued.

We say that the constructible function $\varphi$ on $X$ is {\it
$\pol$-euler} if all functions obtained from $\varphi$ by
using the operations $+$, $-$, $*$, $\hlink$ and $\pol$ are
integer-valued.  The semialgebraic set
$X$ is {\it $\pol$-euler} if $\one_X$ is $\pol$-euler. Thus
a necessary condition for a semialgebraic set $X$ to be homeomorphic
to
an algebraic set is that $X$ is $\pol$-euler.

\subsubsection{A general construction of depth $k$ $\pol$-euler
invariants}\label{Pgen}

In order to produce
a list of $\Z/2$-valued obstructions  the vanishing of which are
necessary and sufficient for $X$ to be $\pol$-euler, we follow
{\it verbatim} subsections \ref{Cons1} and \ref{Cons2}, just replacing the
word
``ring" by ``$\pol$-ring".  Thus let $X\subset \R^n$ be a semialgebraic
set, with $d = \dim X$.
We denote by $\lpol (X)$ the ring all functions obtained from $\one_X$ by
using the operations $+$, $-$, $*$, $\hlink$ and $\pol$.  We
define a sequence of subrings of $\lpol (X)$,
\begin{equation*}
\lpol_0(X)\subset\lpol_1(X)\subset\lpol_2(X)\subset \cdots,
\end{equation*}
where $\lpol_0(X)$ is the ring generated by
$\one_X$, and   for $k\geq0$,
$\lpol_{k+1}(X)$ is the $\pol$-ring generated by $\lpol_k(X)$ and
$\{\hlink\varphi\ |\ \varphi\in\lpol_k(X)\}$.

Suppose that $\lpol_k(X)$ is integer-valued, and
suppose we have subsets $G_0,\dots,G_{k-1}$ of
$\lpol_k(X)$ such that for $j=0,\dots,k$, the set $G_0\cup
\cdots\cup G_j$ additively generates the ring
$\lpol_j(X)/(\ideal (X)\cap\lpol_j(X))$, and for all
$\varphi\in G_j$, $\supp\varphi\subset X_{d-j}$.
 We call the conditions that  $\varphi\in G_k$ are euler
 the {\it depth} $k+1$ {\it $\pol$-euler conditions} for $X$.

The compact semialgebraic set $L$ is the link of
a point in a $\pol$-euler space $X$ if and only if the following
conditions are satisfied:

(A) $L$ is $\pol$-euler,

(B) For all $\varphi\in\lpol(L)$, $\int\varphi\dchi$ is even.

We define the {\it characteristic numbers of depth $k$} by
the same construction as in section \ref {Cons2}.

\subsubsection{Invariants in dimension 4} \label{Pinv}
Let $L$ be a compact semialgebraic set of dimension $\le 3$.  We
work out the conditions (A) and (B) for $L$ to be a link of a
$\pol$-euler space.

(A) {\it $L$ is $\pol$-euler.} It is easy to see that
\begin{equation*}
\lpol(L)/(\ideal (L)\cap\lpol(L)) = \hlink(L)/(\ideal
(L)\cap\hlink(L)),
\end{equation*}
so $L$ is $\pol$-euler if and only of $L$ is completely euler.
Indeed $\lpol_1 (L)/(\ideal(L)\cap\lpol(L))$ is generated additively by
$\one_L$,
$\varphi$, $\varphi^2$, $\varphi^3$, since
$2^{[p/2]}f_p(\varphi)$, $p\ge 4$,
 are divisible by $4$ and therefore in $\ideal (L)$.  Similarly
 $\lpol_2(L)/(\ideal (L)\cap\lpol_2(L))$ is additively
generated by $\one_L$, $\varphi$, $\varphi^2$, $\varphi^3$, and
the products $\alpha\beta$, $\alpha\gamma$,
$\beta\gamma$, $\alpha\beta\gamma$, where
$\beta =\hlink(\varphi^2)$ and $\gamma=\hlink(\varphi^3)$.

(B){\it For all $\varphi\in\lpol(L)$, $\int\varphi\dchi$ is even.}
We compute $\lpol_3(L)/(2\ideal (L)\cap\lpol_3(L))$.
By Theorem \ref{t:pol} (see also \eqref{e:mod8}),
the ring $\lpol_1(L)/(2\ideal (L)\cap\lpol_1(L))$
is additively generated by
\begin{equation}\label{e:pvarphi}
\begin{split}
&\one_L, \,  \varphi, \, \varphi_2 =(\varphi^2 - \varphi),\,
\varphi_3 =(\varphi^3 - \varphi), \\
& \varphi_4 = \half \varphi (\varphi -1)(\varphi-2)(\varphi-3),\,
\varphi_5 = \varphi_4 (\varphi-4) .
\end{split}
\end{equation}
Recall that $\varphi_2\equiv \varphi_3\equiv 0 \pmod 2$ and
$\varphi_4\equiv \varphi_5\equiv 0 \pmod 4$.  So there are no new
 characteristic numbers of depth one.

Now consider
\begin{equation}\label{e:pbcde}
\beta =\hlink(\varphi_2), \,
\gamma =\hlink(\varphi_3),\,
\delta  = \hlink(\varphi_4),\,
\epsilon  = \hlink(\varphi_5).
\end{equation}
(This coincides with our previous notation: $\beta=\hlink(\varphi^2)$ and
$\gamma=\hlink(\varphi^3)$, since $\hlink \varphi =0$.)
As before, let $\alpha = \varphi |L_1$,
$\alpha_1=\alpha$, $\alpha_2=\alpha^2-\alpha$,
$\alpha_3=\alpha^3-\alpha$, and similarly  define
$\beta_1$, $\beta_2$, $\beta_3$, $\gamma_1$, $\gamma_2$,
$\gamma_3$.  Let $\alpha_0= \beta_0 =\gamma_0=1$.
The ring $\lpol_2(L)/(2\ideal (L)\cap\lpol_2(L))$ is
generated
additively by $\one_L$, $\varphi$, $\varphi^2$, $\varphi^3$, together with
representatives of the equivalence classes mod 4 of
\begin{equation}
\psi_{abcde} = \alpha_a\beta_b\gamma_c\delta^d\epsilon^e  \quad b+c+d+e>0 ,
\end{equation}
where $a, b, c=0, 1, 2, 3$, and
$d,e=0,1$ (since $\delta$ and $\epsilon$ are even
valued).  The functions $\psi_{abcde}$ divisible by $4$ do not count, so
the
complete list of generators comprises $S_1\cup S_2\cup S_3\cup
S_4$ of  (\ref{Inv4}) and the following sets:
\begin{equation}
\delta,\, \epsilon;
\tag{$S_0$}
\end{equation}
\begin{gather*}\tag{$S_5$}
\alpha \delta ,\, \beta \delta, \, \gamma\delta,\,
\alpha \epsilon ,\, \beta \epsilon, \, \gamma\epsilon, \\
\alpha\beta\delta, \, \alpha\gamma\delta, \, \beta\gamma\delta, \,
\alpha\beta\epsilon, \, \alpha\gamma\epsilon, \, \beta\gamma\epsilon, \,
  \alpha\beta\gamma\delta,\, \alpha\beta\gamma\epsilon .
\end{gather*}

The characteristic numbers of depth 2 are the euler integrals mod 2
of these functions, and again the only nontrivial ones are
$\int\alpha \beta\dchi ,\, \int\alpha\gamma\dchi ,\,
\int\beta \gamma\dchi,\, \int \alpha \beta \gamma\dchi$.

Let $S = S_0 \cup S_1\cup S_2\cup S_3\cup S_4\cup S_5$, and let
$S'= S_2\cup S_3\cup S_4\cup S_5$.
The ring $\hlink_3(L)/(2\ideal (L)\cap\hlink_3(L))$
is generated additively by $\one_L$, $\varphi$, $\varphi^2$, $\varphi^3$,
and
$S$, together with the equivalence classes mod 2 of the functions
\begin{equation}
  \label{e:pprod}
\Phi(\mathbf m,\mathbf n)=
  \alpha ^{m_\alpha}\beta ^{m_\beta}\gamma ^{m_\gamma}
\prod_{\psi\in S'} (\ahlink (\psi)) ^{n_\psi} ,
\end{equation}
where the exponents $m_\alpha, m_\beta, m_\gamma,n_\psi$ equal $0$ or
$1$, and $\sum_{\psi\in S' }n_\psi >0$.

\begin{rem}\label{r:delta}
Neither $\delta$, $\epsilon$ nor $\ahlink(\delta)$,
$\ahlink(\epsilon)$ appear as factors in \eqref{e:pprod}, since
$\delta$, $\epsilon$ are even valued and
$\ahlink(\delta)= \ahlink(\epsilon)=0$.
\end{rem}

The support of
$\Phi(\mathbf m,\mathbf n)$ is contained in $L_0$,
and hence
\begin{equation} \label{e:pinv}
  \int \Phi(\mathbf m,\mathbf n) \dchi =
\sum_{p\in L_0} \Phi(\mathbf m,\mathbf n)(p) .
\end{equation}
If $\psi\in S_3\cup S_4\cup S_5$ then $\int\psi\dchi$ is
 automatically
even.  This leaves us with
$2^3(2^{40}-1) - 40= 2^{43}-48$ new characteristic numbers
(\ref{e:pinv}) of depth 3.  Thus we obtain:

\begin{thm}\label{t:pinv}
Let $L$ be a compact semialgebraic set of dimension at most $3$.
Let $\varphi = \ahlink \one_L$, $\beta = \hlink \varphi ^2$,
 $\gamma = \hlink \varphi ^3$, $\delta =\hlink \varphi_4$,
$\epsilon = \hlink \varphi_5$.  Then $L$ is
a link in a $\pol$-euler space if and only if the following
conditions hold:
\begin{enumerate}
\item
$L$ is euler and  $\varphi \beta$, $\varphi \gamma$,
$\beta \gamma$,  $\varphi \beta \gamma$
 are euler;
\item
$L$ has even euler characteristic and
   $\varphi \beta$, $\varphi \gamma$, $\beta \gamma$,
  $\varphi \beta \gamma$  
  have even  euler integral;
\item
The $2^{43}-48$ characteristic numbers of (\ref{e:pinv})
\begin{eqnarray*}
\cn (\mathbf m,\mathbf n)(L) =
  \sum_{p\in L_0} \Phi(\mathbf m,\mathbf n)(p),
\end{eqnarray*}
are even.
\end{enumerate}
\end{thm}

\subsection{Independence of invariants}\label{Pind}
We show that the $2^{43} - 43$ characteristic numbers of (2) and (3) of
Theorem \ref{t:pinv} are independent.  For each characteristic
number we construct a
set $L$, together with a filtration $L_0\subset L_1\subset L_2 \subset
L_3=L$ by skeletons of a stratification,
with exactly the given characteristic number nonzero.
For the characteristic number
$\chi (L) \pmod 2$ we can again use $L=S^3\vee S^3$.

\subsubsection{Step 1}\label{Pstep1}
Given a characteristic number $\cn (\mathbf m,\mathbf n)$,
we construct $L_1$ and constructible functions $\alpha$, $\beta$,
$\gamma$, $\delta$, $\epsilon$ on $L_1$ such that:
\begin {enumerate}
        \item[(a)]
        $\ahlink \alpha = \ahlink \beta =\ahlink \gamma =
        \ahlink \delta= \ahlink
        \epsilon =0$, $\alpha \beta$, $\alpha
        \gamma$, $\beta \gamma$, $ \alpha \beta \gamma$, and
        $\one_{L_1}$ are euler,
        and $\delta$, $\epsilon$ are even-valued,
        \item[(b)]
        $\sum_{p\in L_0} \Phi(\mathbf m,\mathbf n)(p)$
        is the only sum as in (\ref {e:pinv}) which is nonzero (mod 2).
\end{enumerate}
We use the elementary blocks $(\alpha)$, $(\beta)$, $(\gamma)$,
$(\alpha\beta),\dots,$ of section \ref{Step1}, with the
understanding that the
values of $\delta$ and $\epsilon$ on these blocks are identically
equal to zero.  We introduce new blocks for the functions in
$S_5$.
(For the functions in $S_0$ we do not need elementary blocks because
of Remark \ref{r:delta}.)
Since $\ahlink \delta = \ahlink \epsilon =0$, the values of
$\delta$ and $\epsilon$ at the vertices are determined by their
values on open 1-dimensional strata.

\begin{rem}\label{r:pcycle}
There is an important difference between $\delta$ and $\epsilon$
and the other even-valued functions of $S$, since the sets
$\{p\in L_1\ |\ \delta/2 \equiv 1 \pmod 2\}$,
$\{p\in L_1\ |\ \epsilon/2 \equiv 1 \pmod 2\}$ are automatically
1-cycles (mod 2).  This follows immediately from the fact that
$\ahlink (\delta/2) = \ahlink (\epsilon/2) =0$. 
\end{rem}

At each vertex of our new elementary blocks we consider, in addition to
$V_1\cup V_2\cup V_3\cup V_4$,  the following
14 new expressions:
\begin{gather*}\tag{$V_5$}
\ahlink  (\alpha \delta) (p) ,\, \ahlink  (\beta \delta) (p), \,
\ahlink  (\gamma\delta) (p),\, \ahlink  (\alpha \epsilon) (p), \,
\ahlink  (\beta \epsilon) (p), \, \ahlink  (\gamma\epsilon) (p), \\
\ahlink  (\alpha\beta\delta) (p), \, \ahlink  (\alpha\gamma\delta) (p), \,
\ahlink  (\beta\gamma\delta) (p), \, \ahlink  (\alpha\beta\epsilon) (p), \,
\ahlink  (\alpha\gamma\epsilon) (p), \, \ahlink  (\beta\gamma\epsilon) (p),
\\
\ahlink  (\alpha\beta\gamma\delta) (p),\, \ahlink
(\alpha\beta\gamma\epsilon) (p) .
\end{gather*}

The elementary blocks corresponding to $V_5$ are defined by simple pictures.
We illustrate 3 cases; the
remaining 11 elementary blocks
are defined by symmetry among the functions $\alpha$, $\beta,$ and
$\gamma$, and by symmetry between $\delta$ and $\epsilon$.

\unitlength = 1cm
\begin{picture}(12,4)
\thicklines
\put(3.1,2){\circle*{.2}}
\put(6,2){\circle*{.2}}
\qbezier (3,2),(4.5,4.4),(6,2)
\qbezier (3,2),(4.5,-.4),(6,2)
\qbezier (3,2),(4.5,2.7),(6,2)
\qbezier (3,2),(4.5,.8),(6,2)
\put(0,2) {\shortstack{$(\alpha\delta)$}}
\put(1.5,2.8) {\shortstack{$(0,0,0,2,0)$}}
\put(3.5,2.4) {\shortstack{$(1,0,0,2,0)$}}
\put(3.5,1.7) {\shortstack{$(1,0,0,0,0)$}}
\put(1.5,1.1) {\shortstack{$(0,0,0,0,0)$}}
\put(6.3,2) {\shortstack {$p$}}
\put(9,2) {\shortstack{$\ahlink(\alpha\delta) (p) \ne 0$}}
\end{picture}

\unitlength = 1cm
\begin{picture}(12,4)
\thicklines
\put(3.1,2){\circle*{.2}}
\put(6,2){\circle*{.2}}
\qbezier (3,2),(4.5,4.4),(6,2)
\qbezier (3,2),(4.5,-.4),(6,2)
\qbezier (3,2),(4.5,2.7),(6,2)
\qbezier (3,2),(4.5,.8),(6,2)
\put(0,2) {\shortstack{$(\alpha\beta\delta)$}}
\put(1.5,2.8) {\shortstack{$(0,0,0,2,0)$}}
\put(3.5,2.4) {\shortstack{$(1,1,0,2,0)$}}
\put(3.5,1.7) {\shortstack{$(1,1,0,0,0)$}}
\put(1.5,1.1) {\shortstack{$(0,0,0,0,0)$}}
\put(6.3,2) {\shortstack {$p$}}
\put(9,2) {\shortstack{$\ahlink(\alpha\beta\delta) (p) \ne 0$}}
\end{picture}

\unitlength = 1cm
\begin{picture}(12,4)
\thicklines
\put(3.1,2){\circle*{.2}}
\put(6,2){\circle*{.2}}
\qbezier (3,2),(4.5,4.4),(6,2)
\qbezier (3,2),(4.5,-.4),(6,2)
\qbezier (3,2),(4.5,2.7),(6,2)
\qbezier (3,2),(4.5,.8),(6,2)
\put(0,2) {\shortstack{$(\alpha\beta\gamma\delta)$}}
\put(1.5,2.8) {\shortstack{$(0,0,0,2,0)$}}
\put(3.5,2.4) {\shortstack{$(1,1,1,2,0)$}}
\put(3.5,1.7) {\shortstack{$(1,1,1,0,0)$}}
\put(1.5,1.1) {\shortstack{$(0,0,0,0,0)$}}
\put(6.3,2) {\shortstack {$p$}}
\put(9,2) {\shortstack{$\ahlink(\alpha\beta\gamma\delta) (p) \ne 0$}}
\end{picture}

Let
$V=V_1\cup V_2\cup V_3\cup V_4 \cup V_5$.  For each
expression $v$ in $V$ the corresponding  elementary block again has
the property that
$v$ is nonzero at all vertices of the block, and all  expressions
which follow $v$ in the linear ordering of $V$
vanish at both vertices of the block.

The rest of step 1 goes exactly as in step 1 of section
\ref{Cons}, since the expressions of $V_5$ are additive under
wedging.

\subsubsection{Step 2}\label{Pstep2}
Given $(L_1, \alpha, \beta, \gamma, \delta, \epsilon)$ satisfying (a) of step
1, we construct $L_2$ and a constructible function $\varphi$ on $L_2$ such
that:
        \begin {enumerate}
        \item[(a)]
        $\hlink \varphi =0$ ,
        \item[(b)]
        $\alpha \equiv \varphi |_{L_1}$, $\beta \equiv \hlink
        \varphi_2$, $\gamma
        \equiv \hlink \varphi_3$, $\delta \equiv \hlink \varphi_4$,
        $\epsilon \equiv \hlink \varphi_5 \pmod {2\ideal (L)}$,
        \item[(c)]
        $\int \varphi \dchi$ is even.
        \end{enumerate}

Let $(L_2, \varphi)$ be such that $\hlink
\varphi =0$, and let $S\subset L_2$ be a 1-dimensional stratum.
Suppose that $S$ is in the boundary
of exactly $k$ two-dimensional strata $S_1, \ldots, S_k$, with the
values of $\varphi$ on the strata equal to $l_1, \ldots,l_k$, respectively.
Let $p\in S$.  Then, since $\hlink \varphi = 0$,
\begin{equation} \label {e:pbetas}
\begin{split}
& \alpha(p) = \varphi (p) = \half \sum l_i, \\
& \beta (p) = \hlink \varphi_2(p) = \varphi_2 (p)
- \half \sum (l^2_i-l_i), \\
& \gamma (p) =  \hlink \varphi_3(p) = \varphi_3 (p)
- \half \sum (l^3_i-l_i) , \\
& \delta (p) = \varphi_4(p) - {\frac 1 4} \sum
l_i(l_i -1)(l_i-2)(l_i-3) , \\
& \epsilon (p) = \varphi_5(p) - {\frac 1 4} \sum
l_i(l_i -1)(l_i-2)(l_i-3)(l_i-4) .
\end{split}
\end{equation}

\begin{lem}\label{l:abcde}
Given
$(a,b,c,d,e)\in(\Z/2)^5$, there exist $l_1,\dots,l_k$,
such that  $(\alpha,\beta,\gamma,\delta,\epsilon)$ given by
\eqref{e:pbetas} satisfy
$(\alpha,\beta,\gamma,\half \delta, \half \epsilon)\equiv
(a,b,c,d,e) \pmod 2$.
\end{lem}

\begin{proof} We have
$\varphi_2 \equiv \varphi_3\equiv \half \varphi_4\equiv \half
\varphi_5\equiv 0 \pmod 2$ and consequently
$\beta$, $\gamma$, $\half \delta$, $\half
\epsilon \pmod 2$ are additive with respect to the operation of
taking the wedge at $p$.  So is $\alpha$.  Let
$P_2(t) = t^2 -t$, $P_3(t) = t^3 -t$, $P_4(t) = \half t(t-1)(t-2)(t-3)$,
$P_5(t) = P_4(t) (t-4)$.  Then
\begin{equation}\label{e:abcde}
\begin{split}
& (\half P_2(1),\half P_3(1), \half P_4(1), \half P_5(1))
\equiv (0,0,0,0) \pmod 2 , \\
& (\half P_2(2),\half P_3(2), \half P_4(2), \half P_5(2))
\equiv (1,1,0,0) \pmod 2 , \\
& (\half P_2(3),\half P_3(3), \half P_4(3), \half P_5(3))
\equiv (1,0,0,0) \pmod 2 , \\
& (\half P_2(4),\half P_3(4), \half P_4(4), \half P_5(4))
\equiv (0,0,1,0) \pmod 2 , \\
& (\half P_2(5),\half P_3(5), \half P_4(5), \half P_5(5))
\equiv (0,0,1,1) \pmod 2 .
\end{split}
\end{equation}
So the lemma follows easily from the additivity property.
\end{proof}

\begin{prop}\label{p:pcycles}
 Let $(L_1,\alpha,\beta,\gamma,\delta,\epsilon)$ be such that
$\ahlink\alpha=
 \ahlink\beta=\ahlink\gamma=0$, $\alpha\beta$, $\alpha\gamma$,
 $\beta\gamma$, $\alpha\beta\gamma$, $\one_{L_1}$ are euler, and
$\delta$, $\epsilon$ are even-valued.
 There exists $(L_2,\varphi)$ such that
\begin{enumerate}
\item[(a)]
$\hlink \varphi = 0$,
\item[(b)]
$\alpha \equiv \varphi |_{L_1}$, $\beta \equiv \hlink
\varphi^2$, $\gamma \equiv \hlink \varphi^3$,
$\half \delta \equiv \half \hlink
\varphi_4$, $\half \epsilon \equiv \half \hlink \varphi_5
\pmod 2$ on open 1-dimensional strata.
\end{enumerate}
\end{prop}

\begin{proof}
By Proposition \ref{p:cycles} there exists $(L_2' ,\varphi)$
such that $\hlink \varphi = 0$, and
$\alpha \equiv \varphi |_{L_1}$, $\beta \equiv \hlink
\varphi^2$, $\gamma \equiv \hlink \varphi^3
\pmod 2$ on open 1-dimensional strata. Moreover, the proofs of
Proposition \ref{p:cycles} and
Lemma \ref{l:abcde} actually show the existence of such
$(L_2' ,\varphi)$ with the additional properties that $\hlink
 \varphi_4 = \hlink  \varphi_5 = 0$, for
$( L_2' ,\varphi)$  can be obtained from $L_1$ by attaching
discs $D_i$ with $\varphi|\operatorname{int}D_i$ equal to
1, 2, or 3 (by (\ref{e:abcde})).

Now by Remark \ref{r:pcycle} the sets
$\{p\in L_1\, |\, \delta/2 \equiv 1 \pmod 2\}$,
$\{p\in L_1\, |\, \epsilon/2 \equiv 1 \pmod 2\}$ are
1-cycles (mod 2).  So by attaching further discs $D_i$ to $L_1$ as
in Proposition \ref{p:cycles}, with $\varphi|\operatorname{int}D_i$
equal to 4 or 5
(by (\ref{e:abcde})), we achieve the desired result.
\end{proof}

Now step 2 is completed just as before. By wedging along 1-strata with
bubbles $(B_{l_i},\varphi,S_{l_i})$ as in Proposition \ref{p:bubbles},
we can adjust the values of $\varphi|_{L_1}$, $\hlink\varphi^2$,
$\hlink\varphi^3\pmod 4$ to agree with $\alpha$, $\beta$, $\gamma$,
respectively, without changing $\hlink\varphi_4$ or
$\hlink\varphi_5$. As a result we have $(L_2,\varphi)$ such that
$\hlink \varphi =0$,
$\alpha \equiv \varphi |_{L_1}$, $\beta \equiv \hlink
\varphi_2$, $\gamma
\equiv \hlink \varphi_3$, $\delta \equiv \hlink \varphi_4$,
$\epsilon \equiv \hlink \varphi_5 \pmod 4$ on open
1-dimensional strata. Then, as in the proof of Corollary \ref{c:step2},
we adjust the values of $\varphi$ on $L_0$ so that $\alpha\equiv
\varphi\pmod 2$ on $L_0$, and it follows that the congruences for
$\beta$, $\gamma$, $\delta$, and $\epsilon$ hold mod 2 on $L_0$.
Finally we can achieve that $\int\varphi\dchi\equiv 0\pmod 2$ by
wedging $L_2$ with a 2-sphere at a point.

\subsubsection{Step 3}\label{Pstep3} Given $(L_2,\varphi)$ satisfying
(a) of step 2, we construct $L$ such that $L$ is euler and
$\varphi\equiv\ahlink\one_L\pmod {2\ideal(L)}$. This is
identical to our previous construction (section \ref{Step3}).

This completes the construction of examples, and hence the proof of
the independence of the invariants (2) and (3) of Theorem \ref{t:pinv}.


\end{document}